\documentclass[preprint,twocolumn,5p]{elsarticle}

\usepackage[utf8]{inputenc}
\usepackage{amsmath}
\usepackage{amsfonts}
\usepackage{amssymb}
\usepackage{amsthm}
\usepackage{booktabs}
\usepackage{graphicx}
\usepackage{hyperref}
\usepackage{import}
\usepackage{mathtools}
\usepackage{microtype}
\usepackage{multirow}
\usepackage{pdfpages}
\usepackage{placeins}
\usepackage{stix2}
\usepackage{tikz}
\usepackage{transparent}
\usepackage{xcolor}
\usepackage[normalem]{ulem} %

\usepackage{cleveref}

\title{A differentiable software suite for accelerated simulation of turbulent flows}

\date{\today}

\begin{document}

\begin{frontmatter}

\author[cwi,eindhoven]{Syver Døving Agdestein\corref{cor1}}
\ead{sda@cwi.nl}
\author[cwi,eindhoven]{Benjamin Sanderse}
\ead{b.sanderse@cwi.nl}

\cortext[cor1]{Corresponding author}

\affiliation[cwi]{
    organization={Scientific Computing Group, Centrum Wiskunde \& Informatica},
    addressline={Science Park 123}, 
    city={Amsterdam},
    postcode={1098 XG}, 
    country={The Netherlands}
}
\affiliation[eindhoven]{
    organization={Centre for Analysis, Scientific Computing and Applications, Eindhoven University of Technology},
    addressline={PO Box 513}, 
    city={Eindhoven},
    postcode={5600 MB}, 
    country={The Netherlands}
}

\begin{abstract}
We present \texttt{IncompressibleNavierStokes.jl}, an open-source Julia package
for solving the incompressible Navier--Stokes equations on staggered Cartesian
grids. The package features matrix-free, hardware-agnostic kernels that are
compiled from a single source for multi-threaded CPU or GPU execution, and
hand-written adjoint kernels for all discrete operators, enabling efficient
reverse-mode automatic differentiation through the entire solver. This
differentiability allows neural network closure models to be trained
\emph{a-posteriori} while embedded in a large-eddy simulation. Memory
optimizations permit double-precision direct numerical simulations at
resolutions up to $840^3$ on a single GPU. The software design, numerical
methods, hardware performance, and integration of neural network closure models
are described, and results for turbulent channel flow are validated against
reference data.

\end{abstract}


\end{frontmatter}

\section{Introduction} \label{sec:introduction}

We present \texttt{IncompressibleNavierStokes.jl}%
~\cite{agdesteinScientificSoftwareSuite2026}
(hereafter referred to as \texttt{INS.jl}),
an open-source software package written in the Julia programming
language~\cite{bezansonJuliaFreshApproach2017} for solving the incompressible
Navier--Stokes equations using second-order staggered finite volume
discretizations on Cartesian
grids~\cite{harlowNumericalCalculationTimeDependent1965}.
The package supports both direct numerical simulation (DNS) and large-eddy
simulation (LES) in two and three dimensions, with various boundary conditions
and non-uniform grids.

The package was initially based on the MATLAB codes INS2D/INS3D%
\footnote{%
    \url{https://github.com/bsanderse/INS2D}
},
where sparse matrices were precomputed and stored for all discrete
operators~\cite{sanderseEnergyconservingDiscretizationMethods2013}.
While conceptually straightforward, storing sparse matrices for all the
operators requires substantial memory, limiting the achievable resolution.
\texttt{INS.jl} replaces these with matrix-free,
hardware-agnostic kernels that are compiled for multi-threaded CPU or GPU from
a single source implementation, drastically reducing memory requirements.

The solver is fully differentiable: hand-written adjoint kernels for all
discrete operators are registered through the \texttt{ChainRules.jl} interface,
enabling reverse-mode automatic differentiation through the entire solver
using \texttt{Zygote.jl}~\cite{innesDontUnrollAdjoint2019}. This makes it possible
to train neural network closure models \emph{a-posteriori}, while embedded
in the LES solver.

The package is designed to be easy to modify and extend.
Adding new physics---such as a temperature field for
Rayleigh--B\'{e}nard convection~\cite{sanderseEnergyconsistentDiscretizationViscous2025},
an Ornstein--Uhlenbeck process for stochastic forcing~\cite{hoekstraReducedSubgridScale2026},
an ``evolve-filter-relax'' scheme~\cite{ivagnesNewDatadrivenEnergystable2026},
or a neural network closure model for LES~\cite{agdesteinDiscretizeFirstFilter2025}---requires
only implementing a new body force or right-hand side term with the
appropriate interface.

This package aims to provide the full computational pipeline---from DNS
data generation through closure model training to LES evaluation---in a single
environment and language.
Julia was chosen for its ability to combine the ease of use of a high-level
language with the performance of compiled
code~\cite{bezansonJuliaFreshApproach2017}, while its multiple-dispatch type
system and composable package ecosystem make it straightforward to integrate
GPU computing, automatic differentiation, and neural network libraries
within a single codebase---without the inter-language interfacing
required by frameworks such as
Relexi~\cite{kurzRelexiScalableOpen2022} and
SmartFlow~\cite{xiaoSmartFlowCFDsolveragnosticDeep2025}, which couple
Python-based machine learning libraries with external CFD solvers
through middleware.

The key features that distinguish \texttt{INS.jl} are:
(i) hand-written adjoint kernels for all discrete operators, enabling
efficient reverse-mode differentiation through the entire solver;
(ii) hardware-agnostic kernels compiled from a single source for
multi-threaded CPU or GPU;
(iii) seamless integration of neural network closure models from
\texttt{Lux.jl} within the differentiable solver;
(iv) support for both single and double precision arithmetic; and
(v) memory optimizations that allow double-precision DNS at resolutions
up to $840^3$ on a single GPU.

This article is organized as follows.
\Cref{S:sec:spatial} presents the numerical methods, including the spatial
discretization on a staggered Cartesian grid and the temporal discretization
with pressure coupling.
\Cref{S:sec:implementation} discusses the software design, including hardware-agnostic
kernels, differentiability, and memory optimization.
\Cref{S:sec:neural} describes the integration of neural network closure
models and the differentiable simulation pipeline.
\Cref{S:sec:development} describes the software development practices,
including version control, documentation generation, testing, continuous
integration, release management, and reproducibility.
\Cref{S:sec:channel-flow} compares numerical results for a turbulent channel
flow simulation with reference data.
\Cref{S:sec:hardware} discusses hardware considerations and floating point arithmetic.

\section{Numerical methods} \label{S:sec:spatial}

This section describes the spatial and temporal discretization of the
incompressible Navier--Stokes equations used in \texttt{INS.jl}.

\subsection{Staggered Cartesian grid}

The $d$-dimensional domain
$\Omega \coloneq [a_1, b_1] \times \dots \times [a_d, b_d]$
is partitioned into $N \coloneq (N_1, \ldots, N_d)$ finite volumes
$\Omega_I \coloneq \Delta^1_{I_1} \times \cdots \times \Delta^d_{I_d}$,
where $I = (I_1, \ldots, I_d)$ is a multi-index with half-integer entries
$I_\alpha \in \{1/2, \, 2 - 1/2, \, \ldots, \, N_\alpha - 1/2\}$.
The volume is a product of 1D intervals
$\Delta^\alpha_i = [x^\alpha_{i - 1/2}, \, x^\alpha_{i + 1/2}]$
for $i \in \{1/2, \, 2 - 1/2, \, \ldots\}$,
with volume center coordinates
$x^\alpha_i = (x^\alpha_{i - 1/2} + x^\alpha_{i + 1/2}) / 2$.

Following the staggered grid approach of
Harlow and Welch~\cite{harlowNumericalCalculationTimeDependent1965},
the pressure $p$ is stored at volume centers $x_I$,
while the velocity component $u^\alpha$ is stored at the faces
$\Gamma^\alpha_I = \Omega_{I - h_\alpha} \cap \Omega_{I + h_\alpha}$,
where $h_\alpha = e_\alpha / 2$ is a half-index shift in the
$\alpha$-direction.
This staggered placement naturally couples the pressure and velocity fields,
avoiding the need for artificial stabilization.

To avoid index offset computations inside the computational kernels, we use
the convention that each finite volume has exactly one pressure component and
one of each velocity component (defined to the right of the volume center in
their canonical direction). If a velocity component is required to the left
of a volume near the boundaries, an additional ghost volume is added.
Boundary values are included in the solution vectors, and before applying an
operator kernel to a field, the boundary values are filled in (for example
using a periodic extension). The memory footprint of including a ghost volume
in a given direction is negligible.

\subsection{Discrete operators}

All discrete operators are built from the discrete derivative in the
$\alpha$-direction,
\begin{equation} \label{S:eq:discrete-derivative}
    (\delta_\alpha \varphi)_I
    = \frac{\varphi_{I + h_\alpha} - \varphi_{I - h_\alpha}}
    {|\Delta^\alpha_{I_\alpha}|},
\end{equation}
which is a second-order approximation of
$\partial \varphi / \partial x^\alpha$.
Using this notation, the four operators appearing in the semi-discrete
equations \eqref{P2:eq:momentum_discrete} are:
\begin{itemize}
    \item \textbf{Divergence:}
        $(\delta_\alpha u^\alpha)_I = 0$
        (mass conservation).
    \item \textbf{Diffusion:}
        $(\delta_\beta \delta_\beta u^\alpha)_I$
        (applied twice for second-order derivatives).
    \item \textbf{Convection:}
        $ (\delta_\beta (
        \eta^\text{half}_\beta u^\alpha 
        \eta^\text{lin}_\alpha u^\beta))_I$,
        using a skew-symmetric form following Verstappen and
        Veldman~\cite{verstappenSymmetrypreservingDiscretizationTurbulent2003}
        that ensures discrete conservation of kinetic energy.
        The product $u^\alpha u^\beta$ at off-diagonal positions requires
        interpolation with weights chosen to preserve skew-symmetry.
    \item \textbf{Pressure gradient:}
        $(\delta_\alpha p)_I$.
\end{itemize}
All operators are divided by the velocity volume sizes, giving them the
same units as their continuous counterparts. The explicit stencil
expressions for each operator, including the interpolation weights for
the convective term and the adjoint kernels used for automatic
differentiation, are given in \cref{S:sec:stencils}.

\subsection{Semi-discrete equations}

Applying the discrete operators to the integral form of the Navier--Stokes
equations yields a system of ordinary differential equations coupled with
an algebraic
constraint~\cite{verstappenSymmetrypreservingDiscretizationTurbulent2003,sanderseEnergyconservingRungeKutta2013}:
\begin{align}
    D u_h &= 0, \label{S:eq:mass-matrix} \\
    \frac{\mathrm{d} u_h}{\mathrm{d} t}
    &= F(u_h) - G p_h, \label{S:eq:momentum-matrix}
\end{align}
where $D$ is the discrete divergence operator,
$G$ is the pressure gradient operator,
$u_h$ and $p_h$ are the discrete velocity and pressure fields, and
$F(u_h)$ collects the convective, diffusive, and forcing contributions.
The momentum equation is only evaluated in the degrees of freedom of the
velocity field. The velocity field also contains boundary values, for which the
momentum equation is not defined.

The key property of the staggered discretization is that the divergence and
pressure gradient operators are each other's transpose (up to a certain
volume-wise scaling in the non-uniform case). This duality is a discrete
analogue of the continuous integration-by-parts identity and is essential for
discrete energy
conservation~\cite{verstappenSymmetrypreservingDiscretizationTurbulent2003}.

\subsection{Boundary conditions}

The software implements four types of boundary conditions, all through ghost
volumes adjacent to the physical domain:
\begin{enumerate}
    \item \textbf{Periodic.}
        The solution is periodic across the domain boundary.
        Ghost volumes duplicate values from the opposite side.
    \item \textbf{Dirichlet.}
        The velocity is prescribed on the boundary, e.g., $u = 0$
        for no-slip walls. Supports spatially and temporally varying
        boundary data $u = u_\mathrm{BC}(x, t)$.
    \item \textbf{Symmetric.}
        The normal velocity vanishes at the boundary,
        while tangential components have zero normal gradient.
        Used to exploit planes of symmetry.
    \item \textbf{Outlet.}
        The stress is prescribed on the boundary, with zero-gradient conditions on the velocity.
\end{enumerate}

\subsection{Non-uniform grids}

\begin{figure}
    \centering
    \includegraphics[width=0.24\columnwidth]{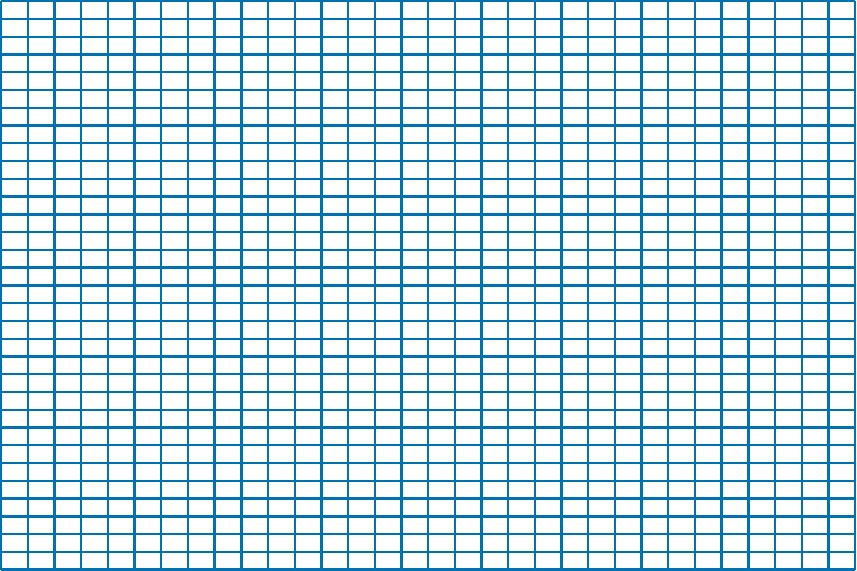}
    \includegraphics[width=0.24\columnwidth]{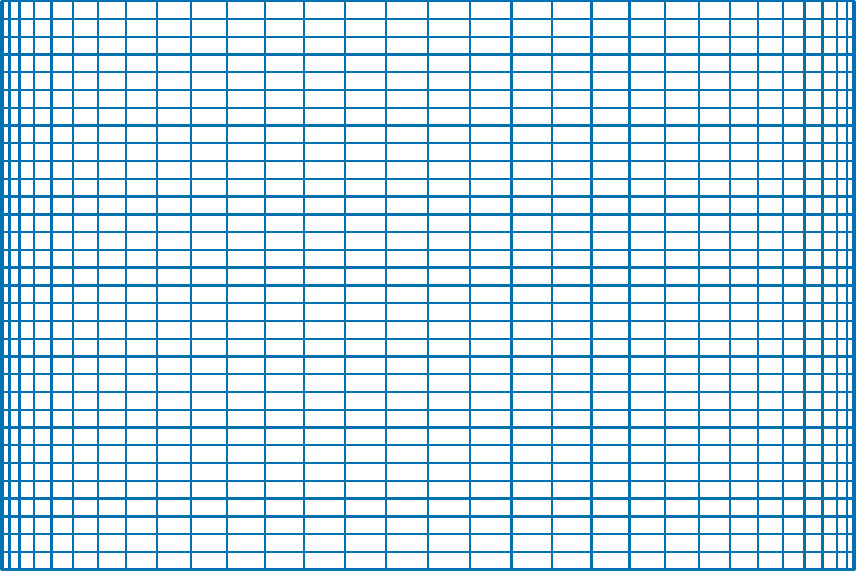}
    \includegraphics[width=0.24\columnwidth]{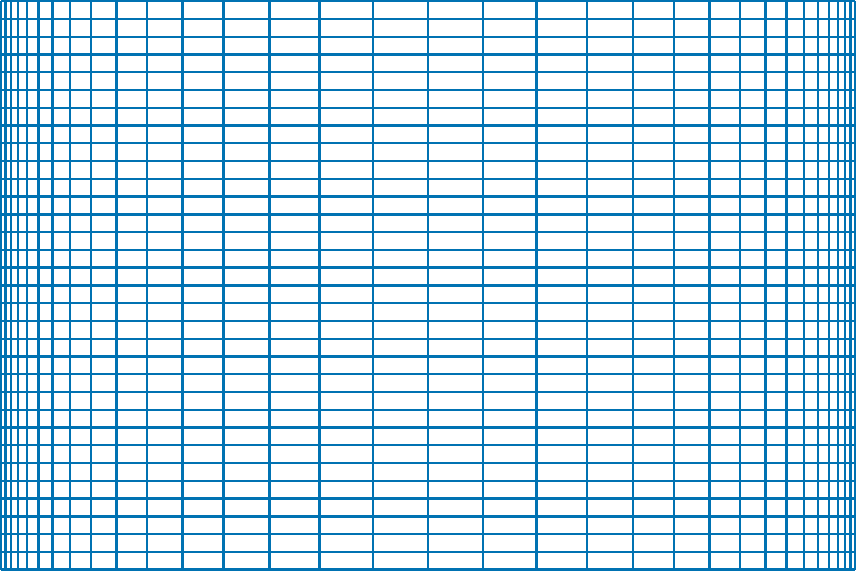}
    \includegraphics[width=0.24\columnwidth]{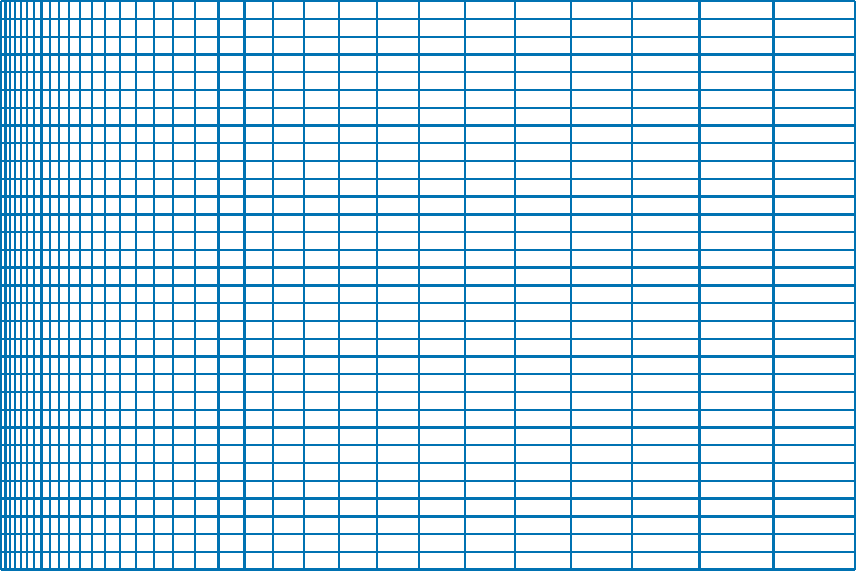}
    \caption{%
        Different horizontal grid spacing profiles. From left to right:
        uniform, cosine, tanh, and stretched grids.
    }
    \label{fig:gridprofiles}
\end{figure}

Non-uniform grids are used to concentrate resolution near boundaries
or regions of interest. Three grid stretching functions are provided (\cref{fig:gridprofiles}).
The \emph{cosine grid} distributes $N + 1$ points on $[a, b]$ as
\begin{equation}
    x_i = a + \frac{1 - \cos(\pi i / N)}{2} (b - a),
    \quad i = 0, \ldots, N,
\end{equation}
which clusters points near both ends.
The \emph{hyperbolic tangent grid}~\cite{triasDirectNumericalSimulations2007}
distributes points as
\begin{equation} \label{S:eq:tanh-grid}
    x_i = a + \frac{b - a}{2}
    \left( 1 + \frac{\tanh\bigl(\gamma (2 \tilde{x}_i - 1)\bigr)}
    {\tanh(\gamma)} \right),
\end{equation}
where $\tilde{x}_i = i / N$ is a uniform parameter and $\gamma > 0$
controls the degree of stretching.
A geometric stretching function
\begin{equation}
    x_i = a + (b - a) \frac{1 - s^i}{1 - s^N}
\end{equation}
is also available for clustering near one end.

\subsection{Temporal discretization and pressure coupling}

Time integration is split into a momentum update and a pressure projection
step to enforce the divergence-free constraint.

\subsubsection{Pressure Poisson equation}

The pressure field is determined at each time step (or sub-step) by
enforcing the divergence-free constraint \eqref{S:eq:mass-matrix}.
Differentiating the mass equation in time and substituting the momentum
equation \eqref{S:eq:momentum-matrix} yields the discrete pressure Poisson
equation
\begin{equation} \label{S:eq:pressure-poisson}
    L p_h = D F(u_h),
\end{equation}
where $L$ is a symmetric positive
(semi-)definite discrete Laplacian. In the absence of pressure boundary
conditions, the pressure is determined only up to a constant, which does
not affect the velocity field.

Three pressure solvers are available:
\begin{enumerate}
    \item \textbf{Spectral solver.}
        For fully periodic domains with uniform grids,
        the Laplace operator $L$ is diagonal in Fourier space.
        The transfer function $\hat{L}$ is written out explicitly, and the
        pressure is computed as
        $p_h = \operatorname{DFT}^{-1} \hat{L}^{-1} \operatorname{DFT} r$,
        where $r$ is the right-hand side and $\operatorname{DFT}$ denotes the
        $d$-dimensional discrete Fourier transform.
        To preserve the average pressure, $\hat{L}_{1, 1}$ is set to one.
    \item \textbf{Direct solver.}
        The Laplace operator is singular, so we directly incorporate the
        zero-mean constraint by solving the augmented system
        \begin{equation}
            \begin{pmatrix} L & e \\ e^\mathsf{T} & 0 \end{pmatrix}
            \begin{pmatrix} p_h \\ \lambda \end{pmatrix}
            = \begin{pmatrix} D F(u_h) \\ 0 \end{pmatrix}.
        \end{equation}
        The augmented matrix is invertible and symmetric, but it is no longer
        positive. We therefore employ the ``square-root free'' Cholesky decomposition
        $L D L^T$ (the $L$ here is not the Laplace matrix, but a lower-triangular factor).
        This requires assembling and storing the sparse Laplace matrix, but
        only a single factorization is needed since the Laplace matrix is
        constant in time.
    \item \textbf{Iterative solver.}
        The conjugate gradient method is applied using the Laplace kernel
        directly (a chain of scaling, divergence, and pressure gradient
        kernels), without assembling the matrix.
        The AMGX solver is also available for GPU.
\end{enumerate}
The Poisson solver is self-adjoint, since $L$ is symmetric. Its pullback
(for reverse-mode automatic differentiation) is therefore itself a Poisson
solve.

\subsubsection{Time integration methods}

The semi-discrete system
\eqref{S:eq:momentum-matrix}--\eqref{S:eq:mass-matrix}
forms a system of differential-algebraic equations.
The software provides a library of time integration methods, including
explicit Runge--Kutta methods (e.g., SSP33, DOPRI6, Wray3) and
the one-leg $\beta$-method of Verstappen and Veldman~\cite{verstappenDirectNumericalSimulation1997}.
At each Runge--Kutta sub-step, the pressure Poisson equation
\eqref{S:eq:pressure-poisson} is solved to project the velocity onto the
divergence-free subspace. This ensures that the velocity is divergence-free
at every sub-step, not only at the end of the time step, and that the
correct order of accuracy of the Runge--Kutta method is obtained.

Adaptive time stepping is supported based on a CFL condition that considers
both convective and diffusive stability limits.

\section{Software implementation} \label{S:sec:implementation}

This section describes the key implementation choices that enable
correctness, performance, and differentiability in \texttt{INS.jl}.

\subsection{Design philosophy}

The discrete operators in the source code map directly to the equations
presented in \cref{S:sec:spatial}. A reader familiar with the mathematical
formulation should be able to read the source code without difficulty.
For example, the divergence operator
$\sum_{\alpha = 1}^d (\delta_\alpha u^\alpha)_I$ is implemented as
the kernel
\begin{verbatim}
divergence(u, D, I) = sum(1:length(I)) do a
    (u[a][I] - u[a][left(I, a, 1)]) / D[a][I[a]]
end
\end{verbatim}
Here, \verb|I| is always a Cartesian index of integers: a quantity
\verb|phi[I]| corresponds to $\varphi_I$ if $\varphi$ is defined at pressure
points, and to $\varphi_{I + h_\alpha}$ if defined at the
$\alpha$-velocity points.

Wherever possible, the code follows a pure functional programming style:
operators are pure functions that take a field as input and return a new
field as output, without modifying any of their arguments. This makes it
straightforward to compose operators, reason about correctness, and reuse
parts of the code in other solvers. Pure functions are also naturally
compatible with reverse-mode automatic differentiation, as discussed in
\cref{S:sec:differentiability}.

However, pure functions allocate a new output array on every call, which
conflicts with the memory constraints of high-resolution DNS
(\cref{S:sec:memory}). To address this, the package also provides
\emph{mutating} variants of all operators, which write their result into
a preallocated output buffer. Following Julia convention, mutating
functions are distinguished by an exclamation mark suffix (e.g.,
\texttt{divergence!} vs.\ \texttt{divergence}). The mutating variants are
used in the time-stepping loop, where memory must be carefully managed,
while the pure variants are used when differentiability is required.
This separation keeps the mutating code isolated and clearly identified,
preserving the benefits of functional style for the majority of the
codebase.

The package is designed for extensibility.
The solver acts on a ``state''---a named tuple of fields that by default
contains only the velocity field. Additional scalar fields (such as a
temperature field for Rayleigh--B\'{e}nard convection) can be added to
this state, and the right-hand side function (excluding the pressure
projection) is user-defined. Adding new physics therefore requires only
implementing a new right-hand side function acting on the augmented
state, without modifying the solver itself. New closure models can
similarly be added by implementing a function with the appropriate
interface.

\subsection{Hardware-agnostic kernels} \label{S:sec:kernels}

All differential operators (divergence, pressure gradient, diffusion,
convection) are implemented as hardware-agnostic kernels using
\texttt{KernelAbstractions.jl}~\cite{churavyJuliaGPUKernelAbstractionsjlV09102023}
as in \texttt{WaterLily.jl}~\cite{weymouthWaterLilyjlDifferentiableBackendagnostic2025}.
The same kernel code is compiled for multi-threaded CPU or GPU, depending on
the type of the input arrays: if called on CPU arrays, the kernel is compiled
to a multi-threaded CPU loop; if called on GPU arrays, it is compiled to a
GPU kernel. The kernels are written in a dimension-agnostic style, so that
the same code handles both two- and three-dimensional problems.

The kernels operate in-place, overwriting preallocated temporary arrays
rather than allocating new memory. These are the mutating operator variants
(with \texttt{!} suffix) described in \cref{S:sec:implementation}; their use is
confined to the time-stepping loop, where memory must be carefully managed.
In-place operation is important for performance, particularly on GPUs where
repeated memory allocation during time stepping would dominate the
computational cost.

\subsection{Memory optimization for single-GPU DNS} \label{S:sec:memory}

A three-dimensional DNS with $N$ grid points per direction requires storing
velocity fields with $d N^d$ components and a pressure field with $N^d$
components, along with temporary arrays for intermediate computations
during time stepping.
The goal is to fit all required arrays in the memory of a consumer GPU,
such as an NVIDIA RTX 4090 with 24 GB of memory, or a datacenter GPU,
such as an NVIDIA H100 with 94 GB of memory.

Three techniques are used to minimize memory consumption:
\begin{enumerate}
    \item \textbf{Array reuse.}
        The mutating operator variants (\cref{S:sec:kernels}) overwrite
        preallocated temporary arrays rather than allocating new memory.
        This is where the departure from pure functional style pays off:
        the same buffer arrays are reused across operators within a single
        time step, avoiding the allocation of new arrays for every
        intermediate result.
    \item \textbf{Low-storage Runge--Kutta methods.}
        Classical Runge--Kutta methods require storing all intermediate
        stages simultaneously, while low-storage variants
        (such as Wray3~\cite{wray1990minimal})
        require only two or three storage registers per degree of freedom,
        independent of the number of stages.
    \item \textbf{Lazy computation of tensor components.}
        Rather than assembling the full $d \times d$ tensor
        $u^\alpha u^\beta$ at once, each component is computed and consumed
        on the fly. This removes the $d^2 N^d$ memory footprint of the
        convective term.
\end{enumerate}
With these techniques, double-precision DNS with $N = 840$ was achieved
on a single H100 in a triply periodic domain.

\subsection{Differentiability and adjoint operators} \label{S:sec:differentiability}

When training neural network closure models, having access to the gradient
of the loss function can greatly speed up the optimization process.
Reverse-mode automatic differentiation (AD) computes these gradients
efficiently by propagating adjoints through the computational graph in
reverse order.

Consider a program $f = f_N \circ \cdots \circ f_1$ composed of elemental
building blocks $f_n$ with intermediate states $x_{n+1} = f_{n+1}(x_n)$.
If $x_N = f(x_1)$ is scalar (e.g., a loss function), the gradient of $f$
with respect to the input $x_1$ can be computed using the pullback functions
\begin{equation} \label{S:eq:pullback-def}
    \bar{f}_{n+1}(x_n) : \bar{x}_{n+1} \mapsto \bar{x}_n
    \coloneq \left( \frac{\mathrm{d} f_{n+1}}{\mathrm{d} x_n}(x_n)
    \right)^\mathsf{T} \bar{x}_{n+1},
\end{equation}
where $\bar{x}_n$ denotes the adjoint (cotangent) variable. Starting from
the seed $\bar{x}_N = 1$, the back-propagation relation
$\bar{x}_n \coloneq \bar{f}_{n+1}(x_n)(\bar{x}_{n+1})$ yields the gradient
$\mathrm{d} f / \mathrm{d} x(x_1) = \bar{x}_1$. Note that the forward states
$x_n$ must be computed first before back-propagation can be performed.

We use \texttt{Zygote.jl}~\cite{innesDontUnrollAdjoint2019} to automatically
generate the back-propagation code. Zygote expands any function $f$ it
encounters until it finds an elemental building block with a
\emph{pullback rule} (or ``rrule'') defined through the \texttt{ChainRules.jl}
framework. Standard operations such as matrix multiplication and neural
network layers have pre-defined rules. However, since our discrete
operators are implemented as matrix-free kernels rather than sparse
matrix multiplications, their pullback rules are not available
automatically. A key contribution of this work is the provision of
hand-written adjoint kernels for all discrete operators, enabling
efficient reverse-mode differentiation through the entire solver.
Because \texttt{Zygote.jl} does not support mutation of arrays, the
pure (non-mutating) operator variants are used during differentiation.
The mutating variants used in the forward time-stepping loop
(\cref{S:sec:kernels}) are bypassed by Zygote in favor of their pure
counterparts, for which the hand-written pullback rules are defined.

Let $\kappa : u \mapsto \varphi$ be a kernel mapping a (possibly vector)
field $u = (u_J)_J$ to a field $\varphi = (\varphi_I)_I$.
The pullback kernel $\bar{\kappa}_u : \bar{\varphi} \mapsto \bar{u}$ at a
given index $J$ is
\begin{equation} \label{S:eq:pullback}
    \bar{u}^\beta_J
    = \left\langle \bar{\varphi},
    \frac{\mathrm{d} \varphi}{\mathrm{d} u^\beta_J}(u) \right\rangle
    = \sum_{\alpha = 1}^{d} \sum_I
    \bar{\varphi}^\alpha_I
    \frac{\mathrm{d} \varphi^\alpha_I}{\mathrm{d} u^\beta_J}(u),
\end{equation}
where $\bar{\varphi}$ is the incoming adjoint and $\bar{u}$ is the
outgoing adjoint.
If $u$ or $\varphi$ is scalar, the corresponding
dimension index and sum are omitted.
Since the kernels are local, most of the derivatives in the pullback are
zero, and the sum over $I$ reduces to a small neighborhood around $J$.

Except for convection, all operators are linear, and their pullbacks reduce
to the transpose of the operator. The convection kernel is nonlinear but
straightforward to differentiate.
The explicit expressions for all kernels and their pullbacks are given in
\cref{S:sec:stencils}.
When a user adds new physics---such as buoyancy, Coriolis terms, or
reaction terms---they can implement the additional terms using standard
vectorized Julia code, which \texttt{Zygote.jl} can differentiate automatically.
Hand-written adjoint kernels are only needed if the user requires the
performance of custom GPU kernels for the new terms; in that case, the
corresponding pullback rules must also be provided.

All adjoint kernels are tested for correctness by comparing against finite
difference gradients. For each operator $\kappa$, we verify that
\begin{equation}
    \langle \bar{\varphi}, \kappa(u + \varepsilon \, \delta u)
    - \kappa(u) \rangle
    \approx \varepsilon \langle \bar{\kappa}_u(\bar{\varphi}), \delta u \rangle
\end{equation}
for random perturbations $\delta u$ and adjoint seeds $\bar{\varphi}$,
at sufficiently small $\varepsilon$.

\section{Neural network closure models} \label{S:sec:neural}

A central motivation for developing \texttt{INS.jl} is the seamless
integration of neural network closure models within the solver,
providing the full computational pipeline---from DNS data generation
through closure model training to LES evaluation---in a single
differentiable framework.

Large-eddy simulation replaces the expensive resolution of all turbulent
scales (DNS) with a coarse-grid simulation supplemented by a closure model
that accounts for the effect of unresolved sub-grid stresses. Classical
closure models based on eddy viscosity, such as the Smagorinsky
model~\cite{smagorinskyGeneralCirculationExperiments1963} and the WALE
model~\cite{nicoudSubgridScaleStressModelling1999}, are implemented in
the software and can serve as baselines---for a full list,
see \cref{S:sec:eddy-viscosity}.

An alternative
to classical closure models is to use neural networks to
learn the closure term from data.
The software supports this through the following components:
\begin{itemize}
    \item Neural network components from \texttt{Lux.jl}, a Julia framework for
        parameterized function approximation.
    \item Reverse-mode automatic differentiation through \texttt{Zygote.jl},
        using the hand-written adjoint kernels described in
        \cref{S:sec:differentiability}.
    \item A-posteriori training: the neural closure model is embedded
        in the LES solver and trained by differentiating through the
        entire time integration, unrolling the solver for a number of
        time steps and backpropagating through them.
\end{itemize}
This differentiable simulation capability distinguishes
\texttt{INS.jl} from most existing CFD solvers, where neural network
training typically requires an external coupling framework. Here, the
hand-written adjoint kernels (\cref{S:sec:differentiability}), the
hardware-agnostic GPU kernels (\cref{S:sec:kernels}), and the neural
network library all operate within the same Julia environment, avoiding
the overhead and complexity of inter-language communication.

A challenge of reverse-mode differentiation through a time-stepping loop
is that \texttt{Zygote.jl} requires the pure (non-mutating) operator variants, which
allocate new arrays for every intermediate result. In the forward
simulation, the mutating variants avoid this overhead
(\cref{S:sec:kernels}), but during training the forward states at each
sub-step must be stored for the backward pass, requiring continuous memory
allocation. The resulting garbage collection overhead can affect
performance, particularly for long training trajectories.

\section{Software development} \label{S:sec:development}

\begin{figure*}
    \centering
    \includegraphics[width=0.8\textwidth]{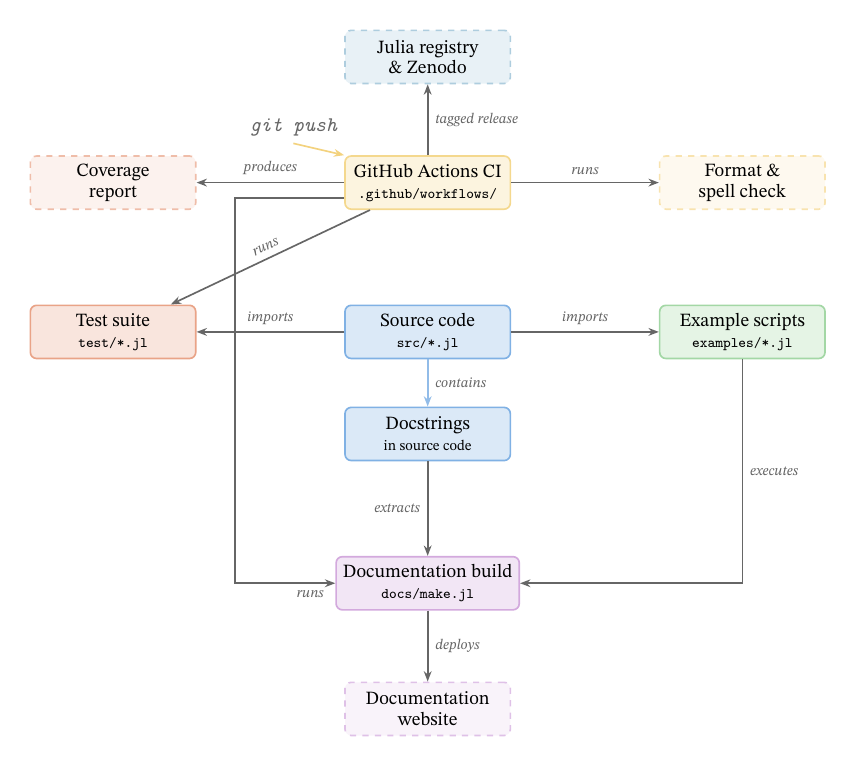}
    \caption{%
        Overview of the software project components and their interactions.
        The source code (center) is imported by the test suite and example
        scripts, and contains docstrings that are extracted during the
        documentation build. Example scripts are executed by \texttt{Literate.jl}
        during the documentation build, producing the documentation website.
        A \texttt{git push} triggers GitHub Actions workflows that run the
        tests, build the documentation, and check code formatting and
        spelling. Tagged releases are published to the Julia General
        registry and archived on Zenodo.
        Solid boxes represent project components; dashed boxes represent
        generated outputs.%
    }
    \label{S:fig:project-structure}
\end{figure*}

While the preceding sections describe the mathematical and computational
content of the software, the present section describes how the software is
\emph{developed}, tested, documented, and released. These practices---version
control, automated testing, continuous integration, and documentation
generation---are standard in professional software engineering but are less
commonly adopted in scientific research, where code is often developed by
individual researchers without formal processes. We describe the practices
used in \texttt{INS.jl} in some detail, as they are essential
for the long-term maintainability and correctness of the code.
An overview of the different components and their interactions is shown in
\cref{S:fig:project-structure}.

\subsection{Version control}

The source code is managed with Git and hosted on GitHub. Every change to
the code is recorded as a commit, providing a complete history of what was
changed, when, and by whom. This makes it possible to trace the origin of
any line of code, identify when a bug was introduced, and revert changes
if needed.

Development follows a branch-based workflow. New features and bug fixes are
developed on dedicated branches (e.g., \texttt{fix-CFL},
\texttt{tensorclosure}) and merged into the main branch via pull requests.
Pull requests allow changes to be reviewed and tested before they are
accepted into the main codebase. Commit messages follow the conventional
commits format (e.g., \texttt{feat: add default time step},
\texttt{fix(examples): minor fixes}, \texttt{docs: update examples to new
syntax}), which makes the history easy to browse and filter by change type.

\subsection{Docstrings}

\begin{figure}
    \centering
    \includegraphics[width=\columnwidth]{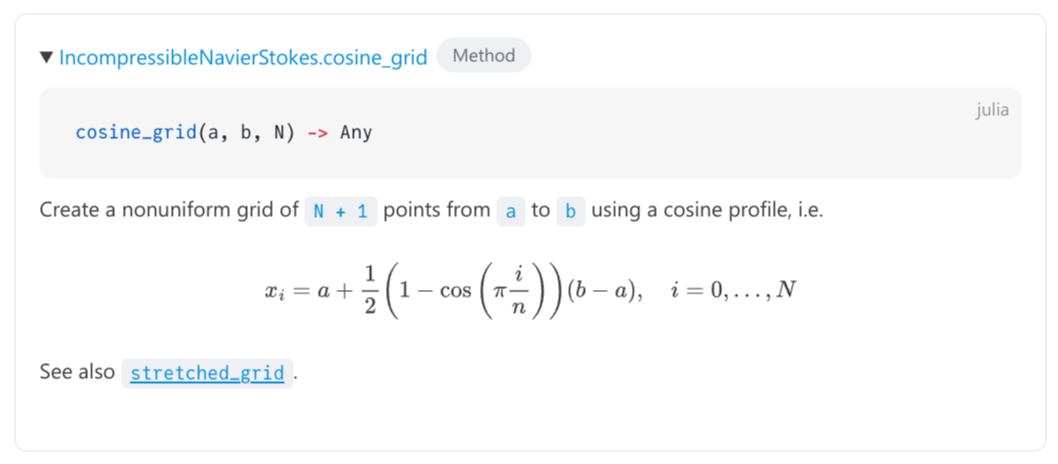}
    \caption{Rendered documentation from a docstring.}
    \label{S:fig:docstring}
\end{figure}

All public functions and types are documented with docstrings: inline
documentation strings placed directly in the source code, adjacent to the
function they describe. In Julia, docstrings are written in Markdown and
placed immediately before the function definition. For example:
\begin{verbatim}
"""
Create a nonuniform grid of `N + 1` points
from `a` to `b` using a cosine profile.
"""
function cosine_grid(a, b, N)
    i = 0:N
    @. a + (b - a) * (1 - cospi(i / N)) / 2
end
\end{verbatim}
Docstrings serve two purposes: they document the code for developers reading
the source, and they are extracted automatically by the documentation
generator to produce a searchable API reference.
This is shown in \cref{S:fig:docstring}.
Cross-references between functions are created using the \texttt{@ref} macro.

\subsection{Documentation generation} \label{S:sec:documentation}

\begin{figure*}
    \centering
    \includegraphics[width=0.8\textwidth]{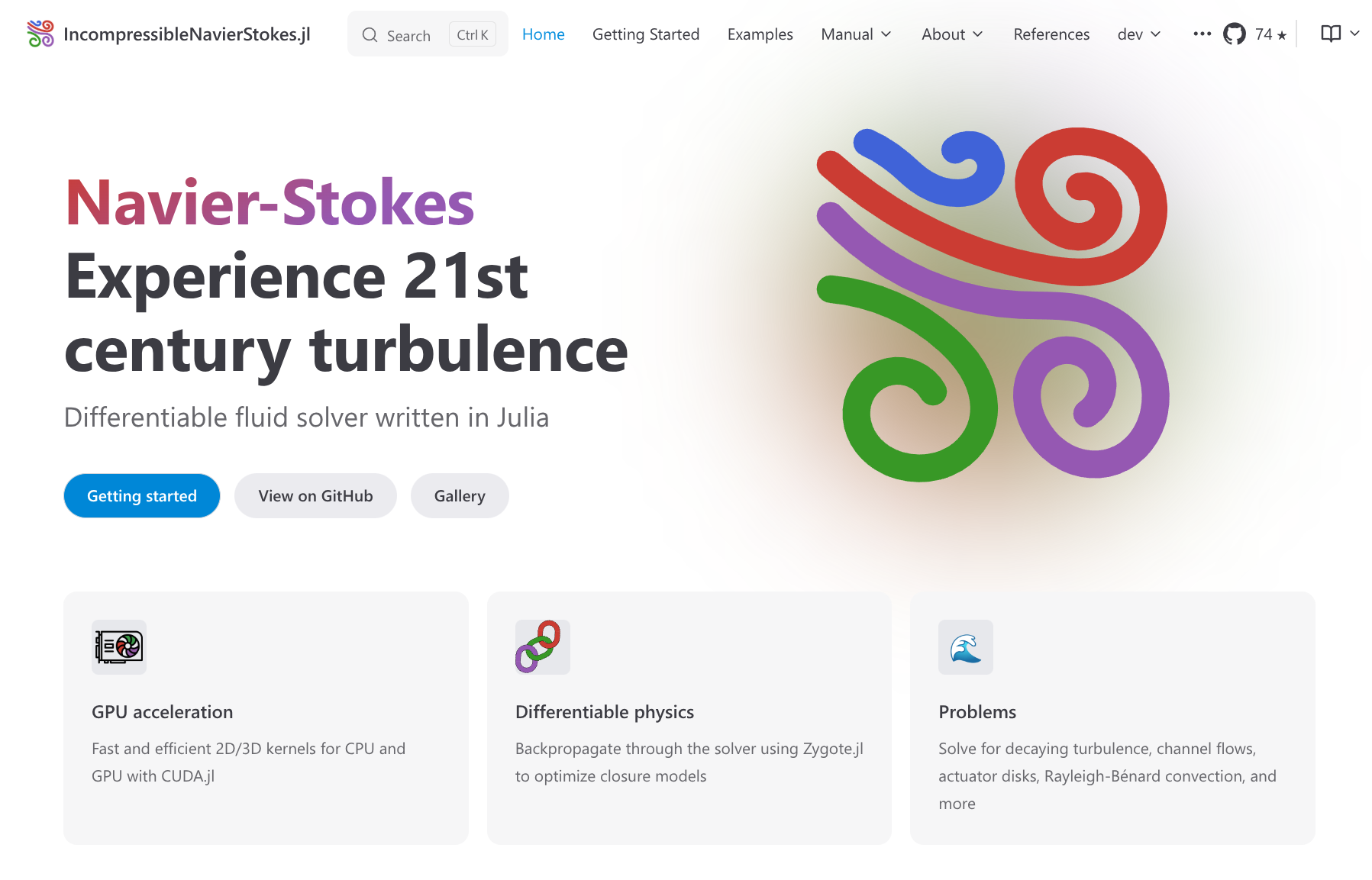}
    \caption{Generated documentation website.}
    \label{S:fig:documentation}
\end{figure*}

The package documentation is generated using
\texttt{Documenter.jl}, which
extracts docstrings from the source code, combines them with hand-written
guide pages, and produces a complete documentation website
(see \cref{S:fig:documentation}).
The website is rendered using \texttt{VitePress} via the \texttt{DocumenterVitepress.jl} backend, and
bibliographic references are handled by \texttt{DocumenterCitations.jl}.

The documentation build is defined in a script (\texttt{docs/make.jl})
that specifies the structure of the documentation, which modules to
extract docstrings from, and how to process example scripts. The
generated website includes an API reference (automatically collected from
docstrings), tutorial examples (generated from literate scripts, described
below), and guide pages.

\subsection{Literate programming for examples} \label{S:sec:literate}

Tutorial examples are written as Julia scripts with embedded Markdown
comments, using \texttt{Literate.jl}.
Lines beginning with \texttt{\#~} are treated as Markdown text, while
the remaining lines are Julia code. During the documentation build,
\texttt{Literate.jl} converts each script into a Markdown page with executable
code blocks. The code is then run by \texttt{Documenter.jl}, and any output---including
plots and figures---is captured and included in the published documentation.
Figures in the examples are produced using
\texttt{Makie.jl}~\cite{danischMakiejlFlexibleHighperformance2021}.
This ensures that all tutorial examples are always up to date with the
current version of the code: if a breaking change is made to the API,
the documentation build will fail, immediately signalling that the
examples need to be updated.

\begin{figure*}
    \centering
    \includegraphics[width=0.7\textwidth]{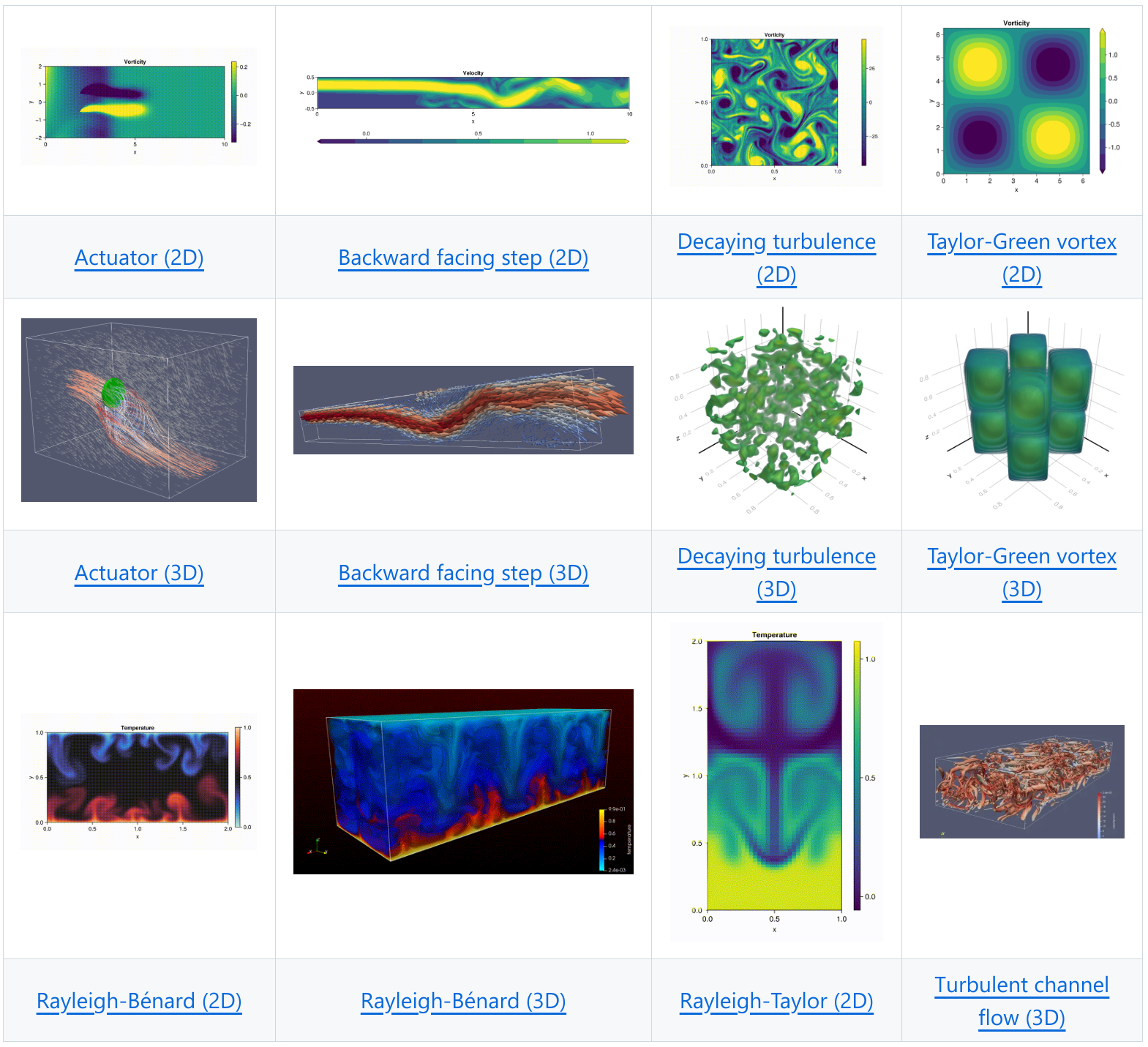}
    \caption{Different example scripts included in the documentation with \texttt{Literate.jl}.}
    \label{S:fig:examples}
\end{figure*}

The package includes fifteen example scripts covering a range of
configurations: periodic turbulence in two and three dimensions,
lid-driven cavity flows, turbulent channel flow,
backward-facing steps, actuator disk models,
Rayleigh--B\'{e}nard convection, and Rayleigh--Taylor instability.
They are shown in \cref{S:fig:examples}.
To manage computational cost, only a subset of the examples is
executed during the documentation build; the remaining scripts are
included as non-executed code listings.

\subsection{Continuous integration} \label{S:sec:ci}

Every push to the main branch, every pull request, and every tagged release
triggers a set of automated workflows on GitHub Actions.
Continuous integration ensures that the code is always in a working state:
any change that breaks the tests, the documentation, or the code formatting
is detected immediately, before it can be merged.

The following workflows are configured:
\begin{enumerate}
    \item \textbf{Unit tests.}
        The full test suite is run on the latest stable Julia version
        on Ubuntu. The workflow has a timeout of 60 minutes to avoid
        runaway computations. Code coverage is measured and uploaded to
        Codecov, to see which lines of code covered by the tests.
        For pull requests, intermediate builds are cancelled when
        a new commit is pushed, avoiding unnecessary resource usage.
    \item \textbf{Documentation.}
        The documentation is built from source, executing the literate
        example scripts (\cref{S:sec:literate}) and extracting docstrings.
        The result is deployed to GitHub Pages. Preview builds are
        generated for pull requests, allowing documentation changes to be
        reviewed before merging.
    \item \textbf{Format check.}
        All source files are reformatted using \texttt{JuliaFormatter.jl}, and the
        workflow fails if any file was changed. This enforces a consistent
        code style across all contributors without requiring manual formatting.
    \item \textbf{Spell check.}
        On pull requests, a spell checker (\texttt{typos}) scans the
        codebase for common typographical errors in code and documentation.
    \item \textbf{Dependency compatibility.}
        \texttt{CompatHelper.jl} runs daily to check whether new versions of
        dependencies are available that may require updating the
        compatibility bounds in \texttt{Project.toml}. When updates are
        needed, it automatically opens a pull request. Similarly,
        Dependabot monitors the GitHub Actions workflow files for outdated
        action versions.
\end{enumerate}
This automated pipeline means that a developer can push a commit and
receive feedback within minutes on whether the change is correct, well
formatted, and compatible with the existing documentation and
dependencies.

\subsection{Release management} \label{S:sec:releases}

The package follows semantic versioning: the version number
(e.g., 3.0.0) encodes major, minor, and patch levels. Major version
increments signal breaking API changes, minor increments add new features,
and patch increments fix bugs.
The version is specified in \texttt{Project.toml}, the central package
manifest.

Releases are published through the Julia General registry. When a new
version is ready, it is registered with the registry, after which
TagBot---a GitHub Actions workflow---automatically creates a GitHub
release and a corresponding Git tag. This makes every release
permanently accessible and installable by any Julia user through the
standard package manager.

Each release is also archived on Zenodo, which assigns a digital object
identifier (DOI) to the archived snapshot. Unlike a GitHub repository,
which can be deleted or made private, a Zenodo archive is intended to be
permanent. This is important for reproducibility: a DOI in a publication
points to the exact version of the code used to produce the results,
regardless of future changes to the repository.

The package also includes a \texttt{CITATION.cff} file, a
machine-readable citation metadata file that specifies the authors,
title, license, and DOI. GitHub recognizes this file and displays a
``Cite this repository'' button on the repository page, making it easy
for users to generate a citation in various formats.

\subsection{Reproducibility} \label{S:sec:reproducibility}

All simulation results presented in scientific computing publications should
be reproducible.
The code used to generate the results should ideally be made publicly
available. If this is not the case, the authors should be willing to help
readers upon demand, or agree to share the code privately upon request with
a non-disclosure agreement.

Simulation setups should be provided in full detail.
While providing a Taylor-scale Reynolds number is descriptive for what kind
of flow is being simulated, it is not sufficient for reproduction, since it
is a \emph{derived} quantity. The simulation setup should also include the
exact domain size, kinematic viscosity, and initial conditions so that
readers can reproduce the simulation setup and compare their results.

When using random numbers, the seed of the random number generators should
always be fixed. This in itself is not a guarantee for perfect
reproducibility, since different hardware architectures and compiler
optimizations may lead to different results. Therefore, the exact hardware
and software environment used to generate results should be explicitly
stated.

Research code used to produce results for publications should be archived
for long-term availability. Zenodo is such an archival service. A GitHub
repository can disappear after some years, but a Zenodo archive is at least
intended to be permanent. Zenodo provides convenient integration with GitHub
for archiving repositories, and each version of the archived software gets
a digital object identifier (DOI).

In addition to archiving code, one may also store the simulation results
themselves. Full three-dimensional velocity fields from DNS or LES require
substantial storage space, which makes long-term hosting nontrivial.
Zenodo offers a solution for this as well, allowing researchers to upload
large datasets with a persistent DOI \cite{hoekstraReducedSubgridScale2026}.
The tradeoff is that raw field data is expensive to store and transfer,
while derived statistics---such as mean velocity profiles, Reynolds stresses,
and energy spectra---are compact and easily shared.
Vreman and Kuerten~\cite{vremanComparisonDirectNumerical2014} stored averaged
DNS statistics for turbulent channel flow in small text files that are readily
downloadable, providing a convenient reference dataset. We will use this
reference data in \cref{S:sec:channel-flow}.

\subsection{Test suite}

The test suite of \texttt{INS.jl} comprises fourteen
test files organized using \texttt{TestItemRunner.jl}. The entire suite is run with
\begin{verbatim}
julia --project=test test/runtests.jl
\end{verbatim}
or equivalently through Julia's package manager with \texttt{Pkg.test()}.
Individual tests can be filtered by name or file, which is useful during
development when working on a specific component.
The tests fall into three categories.

\begin{enumerate}
    \item \textbf{Software correctness tests.}
        These verify that individual functions and complete simulation
        pipelines run without errors. A ``master run'' test steps through
        an entire simulation at low resolution to catch interface-level
        bugs. The numerical values produced are not checked; the purpose is
        to ensure that the code executes without crashing for representative
        use cases.
        Code quality is additionally checked by \texttt{Aqua.jl}, which detects
        common issues such as ambiguous method definitions, unbound type
        parameters, and stale dependencies.
    \item \textbf{Numerical verification tests.}
        These check that the code produces correct numerical values.
        The 2D Taylor-Green vortex serves as a
        manufactured solution: the solver must converge at second order
        toward the known analytical velocity and pressure fields.
        The pressure Poisson solvers are verified by checking that the
        resulting velocity field is divergence-free up to a prescribed
        tolerance.
        All hand-written adjoint kernels (\cref{S:sec:differentiability})
        are compared against finite-difference gradients using
        \texttt{ChainRulesTestUtils.jl}.
    \item \textbf{Physical verification tests.}
        These test that the discrete operators satisfy the physical
        invariants they are designed to preserve.
        The skew-symmetric convection operator must produce zero change in
        the total kinetic energy.
        The discrete divergence and pressure gradient operators must be each
        other's negative transpose.
\end{enumerate}
Where possible, tests are run on both uniform and non-uniform grids and
at small problem sizes, since edge cases---such as very non-uniform grid
spacings or few degrees of freedom---are more likely to expose bugs that
remain hidden on uniform grids at moderate resolution.

All tests are designed to complete within the resource limits of GitHub
Actions (\cref{S:sec:ci}), so that every commit receives automated
feedback within minutes. This keeps the development cycle fast and
ensures that new contributions do not break existing functionality.

Next, we compare the code against a turbulent channel flow to verify
that the numerical solutions produce expected statistical properties.
This is a visual confirmation, and is separate from the automated test
suite.

\section{Turbulent channel flow} \label{S:sec:channel-flow}

\begin{figure}
    \centering
    \includegraphics[width=\columnwidth]{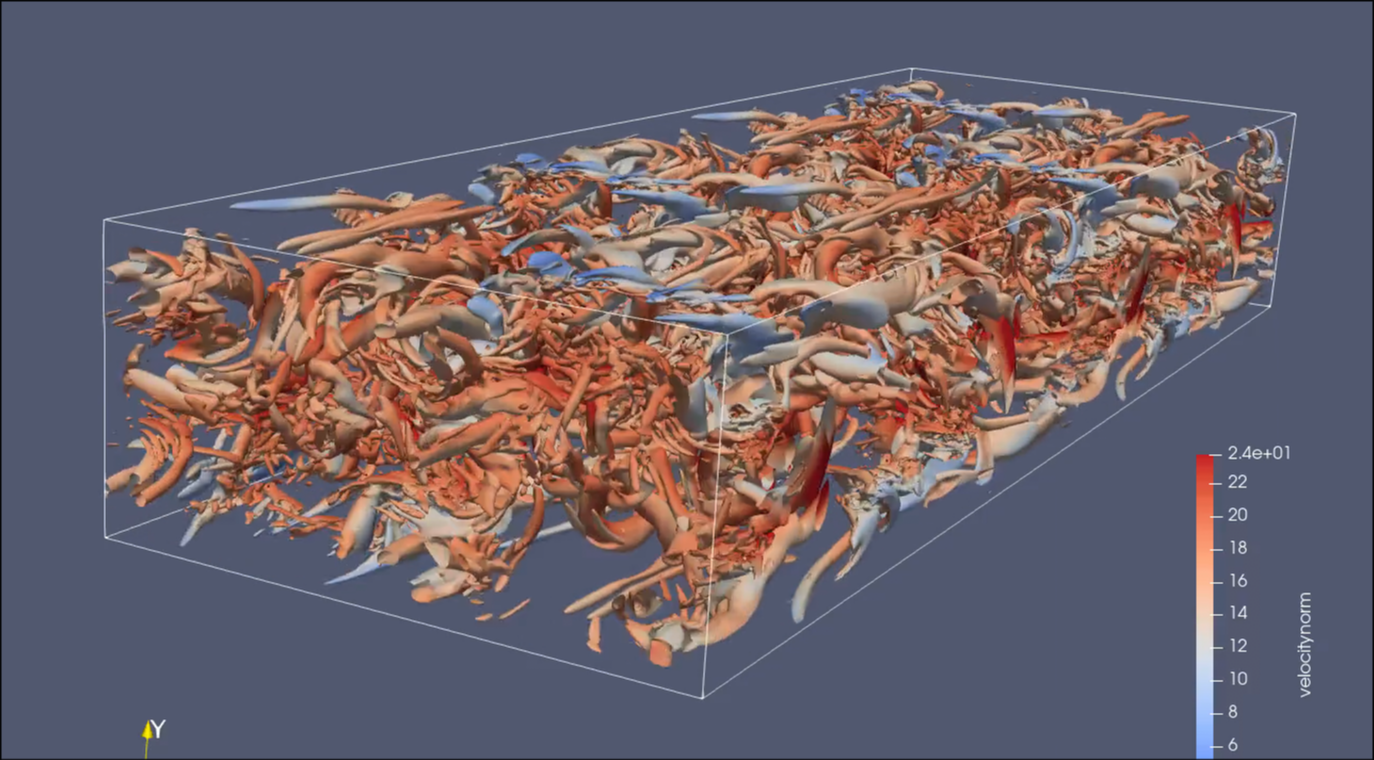}
    \caption{%
        Snapshot of the turbulent channel flow during the transition from
        laminar to turbulent flow. The isosurfaces show the $Q$-criterion
        $Q = -\frac{1}{2} \operatorname{tr}(A^2)$,
        where $A$ is the velocity gradient tensor.
        Positive values of $Q$ identify regions where rotation dominates
        strain, highlighting vortical structures.
        The isosurfaces are colored by velocity magnitude.%
    }
    \label{S:fig:channel-snapshot}
\end{figure}

We simulate turbulent channel flow at friction
Reynolds number $\mathrm{Re}_\tau = 180$ (\cref{S:fig:channel-snapshot})
and compare against the reference DNS data of Vreman and
Kuerten~\cite{vremanComparisonDirectNumerical2014}.

\subsection{Setup}

The domain is a rectangular box
$[0, L_x] \times [0, L_y] \times [0, L_z]$ with
$L_x = 4 \pi$,
$L_y = 2$, and
$L_z = 4 \pi / 3$.
The velocity is denoted $(u, v, w)$.
The boundary conditions are periodic in the streamwise ($x$) and
spanwise ($z$) directions, and no-slip at the walls $y = 0$ and $y = 2$.
The flow is driven by a constant body force $f = (1, 0, 0)$ that acts
as an artificial mean pressure gradient.
The kinematic viscosity is $\nu = 1 / 180$.

The friction velocity is defined as
$u_\tau \coloneq \lim_{y \to 0} \left( \nu \mathrm{d} \bar{u} / \mathrm{d} y \right)^{1/2}$
at the wall,
where $\bar{u}(y) \coloneq \int u(x, y, z) \, \mathrm{d} x \, \mathrm{d} z /
(L_x L_z) $ is the mean streamwise velocity.
In the limit where $T \to \infty$, this setup gives
$\langle u_\tau \rangle = 1$, where $T > 0$ is the averaging time and
$\langle \cdot \rangle$ denotes time averaging.
We therefore use the approximation $u_\tau \approx 1$ to define a velocity scale,
corresponding to the
``wall unit'' length scale $L^+ \coloneq \nu / u_\tau = 1 / 180$,
time scale $T_\tau \coloneq H / u_\tau = 1$, and
Reynolds number
$\mathrm{Re}_\tau \coloneq u_\tau H / \nu = 180$,
where $H = 1$ is the channel half-width.
The initial conditions are a sinusoidally perturbed laminar flow.
The flow is first integrated until $t = 15 T_\tau = 15$ time units to reach
a statistically stationary state, after which statistics are collected at
$1000$ equispaced times for $10 T_\tau = 10$ time units.

\subsection{DNS on a uniform grid}

We first perform a DNS on a uniform grid with
$N_x = 512$, $N_y = 1024$, $N_z = 256$ points.
On this uniform grid, the pressure Poisson equation is solved spectrally
using FFTs in the periodic directions and a discrete cosine transform
(type-II) in the wall-normal direction, implemented via the
\texttt{cuFFT} library on the GPU.
In wall units, this gives grid spacings of
$\Delta x^+ \coloneq \Delta x / L^+ = 4.418$,
$\Delta y^+ = 0.3516$, and
$\Delta z^+ = 2.945$.
Vreman and Kuerten used $N_x = 512$, $N_y = 256$, $N_z = 256$ points with a
$\tanh$-stretching in the wall-normal direction, giving $\Delta y^+ = 0.49$
near the wall
(with the same $\Delta x^+$ and $\Delta z^+$ as our uniform grid).

\begin{figure}
    \centering
    \includegraphics[width=\columnwidth]{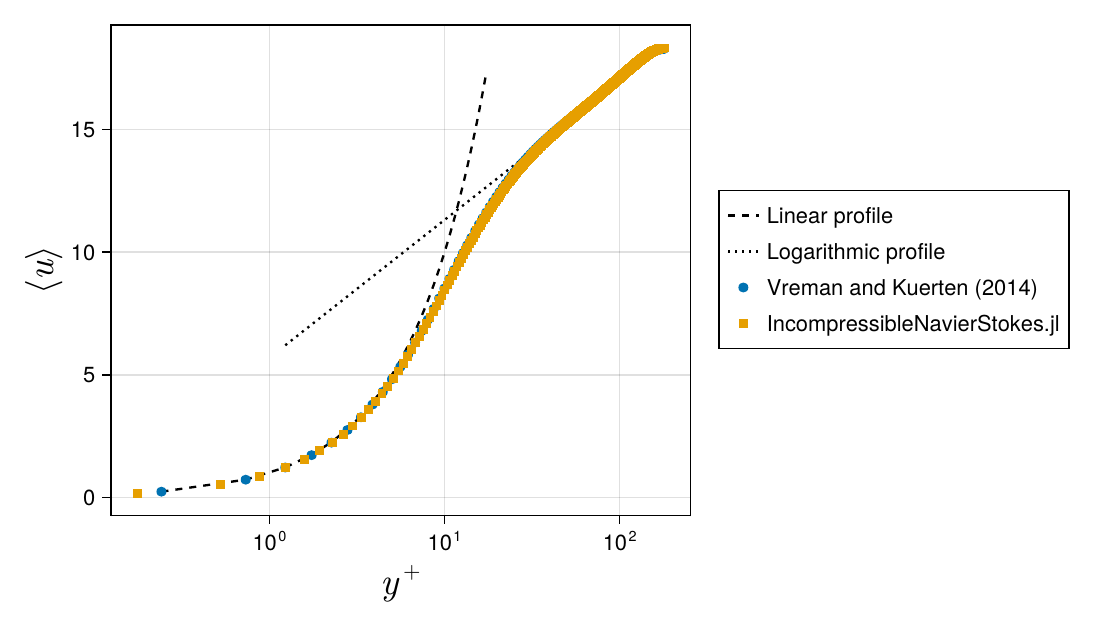}
    \caption{%
        Mean streamwise velocity profile $\langle \bar{u} \rangle$
        for DNS of turbulent channel flow
        with $512 \times 1024 \times 256$ grid points on a uniform
        grid, compared with reference data from Vreman and
        Kuerten~\cite{vremanComparisonDirectNumerical2014}.%
    }
    \label{S:fig:wallplot-dns-ubar}
\end{figure}

\Cref{S:fig:wallplot-dns-ubar} shows the mean streamwise velocity profile
$\langle \bar{u} \rangle$ as a function of the wall-normal coordinate in
viscous units, $y^+ = y / L^+ = 180 \, y$.
The reference data from Vreman and
Kuerten~\cite{vremanComparisonDirectNumerical2014} is also shown.
The velocity profiles are in good agreement
in both the linear viscous sublayer and the logarithmic region.
The theoretical linear ($u / u_\tau = y^+$) and logarithmic
($u / u_\tau = \frac{1}{\kappa} \log y^+ + C$ with $\kappa = 0.41$, $C = 5.7$)
profiles are also shown for reference.

\begin{figure}
    \centering
    \includegraphics[width=\columnwidth]{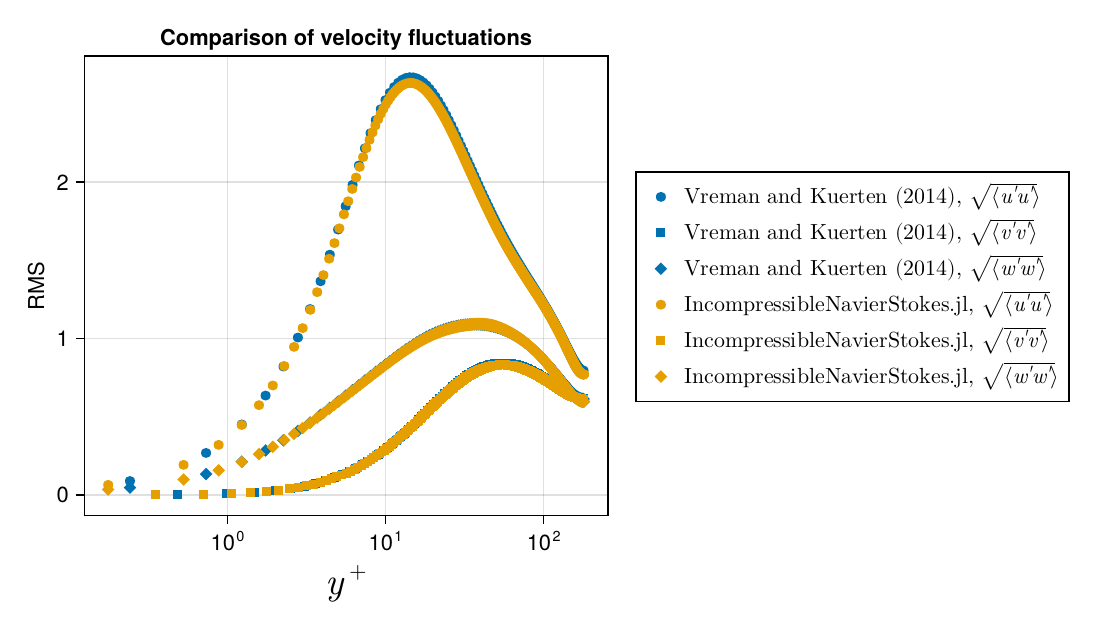}
    \caption{%
        Average velocity fluctuations in the three directions
        for DNS of turbulent channel flow
        with $512 \times 1024 \times 256$ grid points on a uniform
        grid, compared with reference data from Vreman and
        Kuerten~\cite{vremanComparisonDirectNumerical2014}.%
    }
    \label{S:fig:wallplot-dns-rms}
\end{figure}

\begin{figure}
    \centering
    \includegraphics[width=\columnwidth]{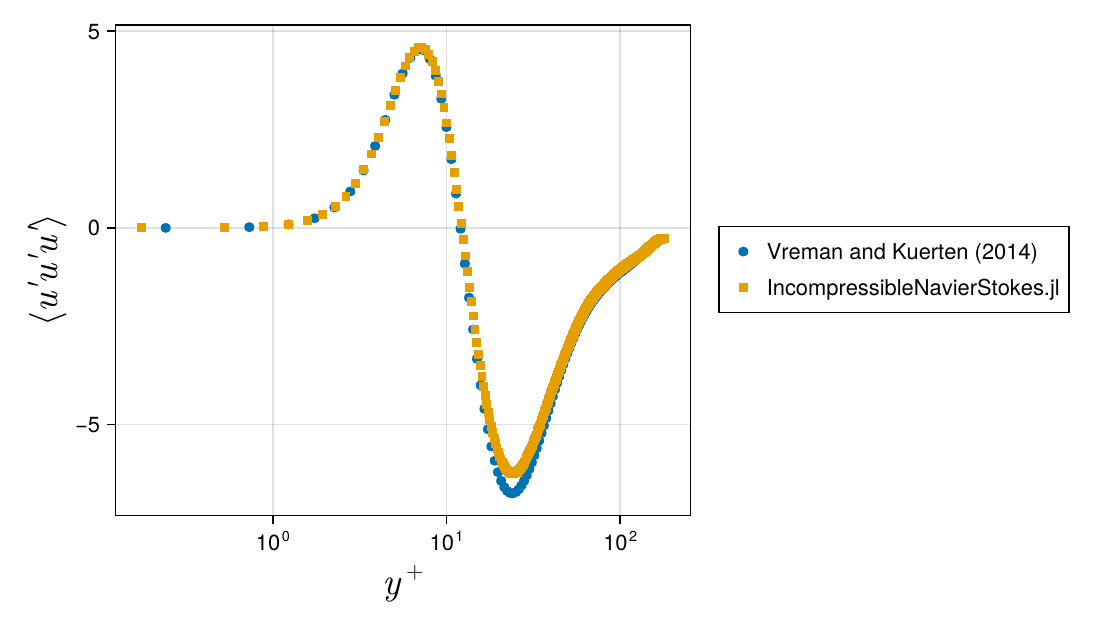}
    \caption{%
        Third moment of $x$-velocity fluctuations $\left\langle \overline{u' u' u'} \right\rangle$
        for DNS of turbulent channel flow
        with $512 \times 1024 \times 256$ grid points on a uniform
        grid, compared with reference data from Vreman and
        Kuerten~\cite{vremanComparisonDirectNumerical2014}.%
    }
    \label{S:fig:wallplot-dns-up3}
\end{figure}

\begin{figure}
    \centering
    \includegraphics[width=\columnwidth]{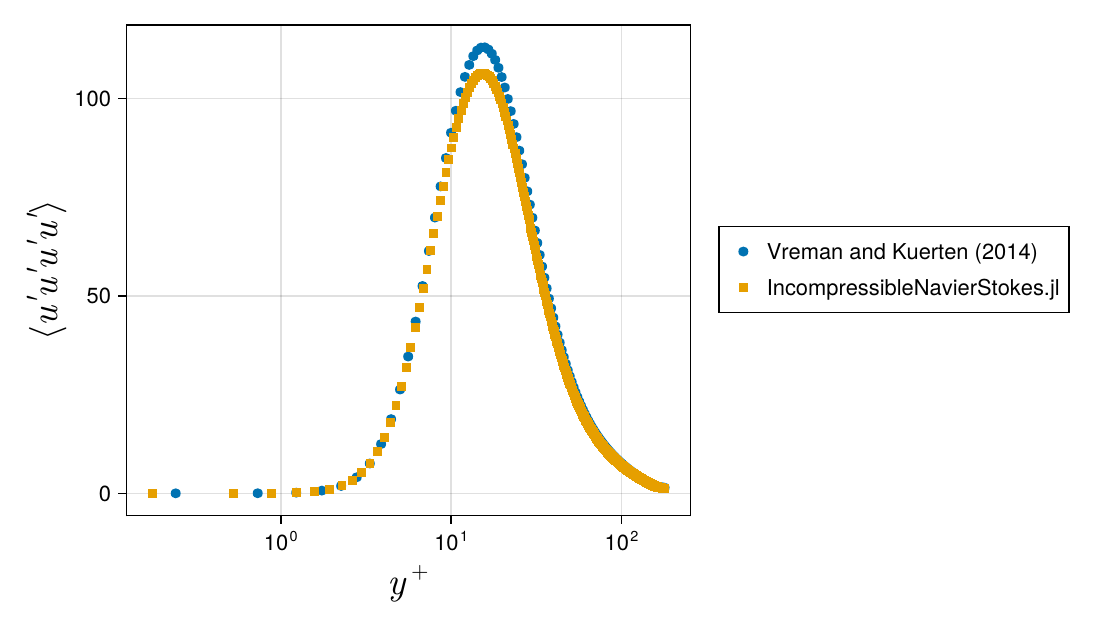}
    \caption{%
        Fourth moment of $x$-velocity fluctuations $\left\langle \overline{u' u' u' u'} \right\rangle$
        for DNS of turbulent channel flow
        with $512 \times 1024 \times 256$ grid points on a uniform
        grid, compared with reference data from Vreman and
        Kuerten~\cite{vremanComparisonDirectNumerical2014}.%
    }
    \label{S:fig:wallplot-dns-up4}
\end{figure}

Further insight can be gained by looking at the higher order moments of the velocity fluctuations
\begin{equation}
    u' \coloneq u - \bar{u}.
\end{equation}
These include the root-mean-square (RMS) velocity fluctuations
$\sqrt{\left\langle \overline{u' u'} \right\rangle}$,
which are shown in \cref{S:fig:wallplot-dns-rms}.
\texttt{INS.jl} is in good agreement with the reference
data, but there is a tiny discrepancy at the peak of
$\left\langle \overline{u' u'} \right\rangle$.
This discrepancy is more visible in the third and fourth moments
$\left\langle \overline{u' u' u'} \right\rangle$ and
$\left\langle \overline{u' u' u' u'} \right\rangle$
shown in \cref{S:fig:wallplot-dns-up3} and
\cref{S:fig:wallplot-dns-up4}, respectively.
The peaks are slightly underpredicted.
We stress that we use a second-order accurate finite-volume scheme, while the
reference data was produced using a fourth-order finite-difference scheme.
Vreman and Kuerten also averaged over $200 T_\tau = 200$ time units, while
we only average over $10$ time units, but we found little difference when
running simulations for longer times.

\begin{figure}
    \centering
    \includegraphics[width=\columnwidth]{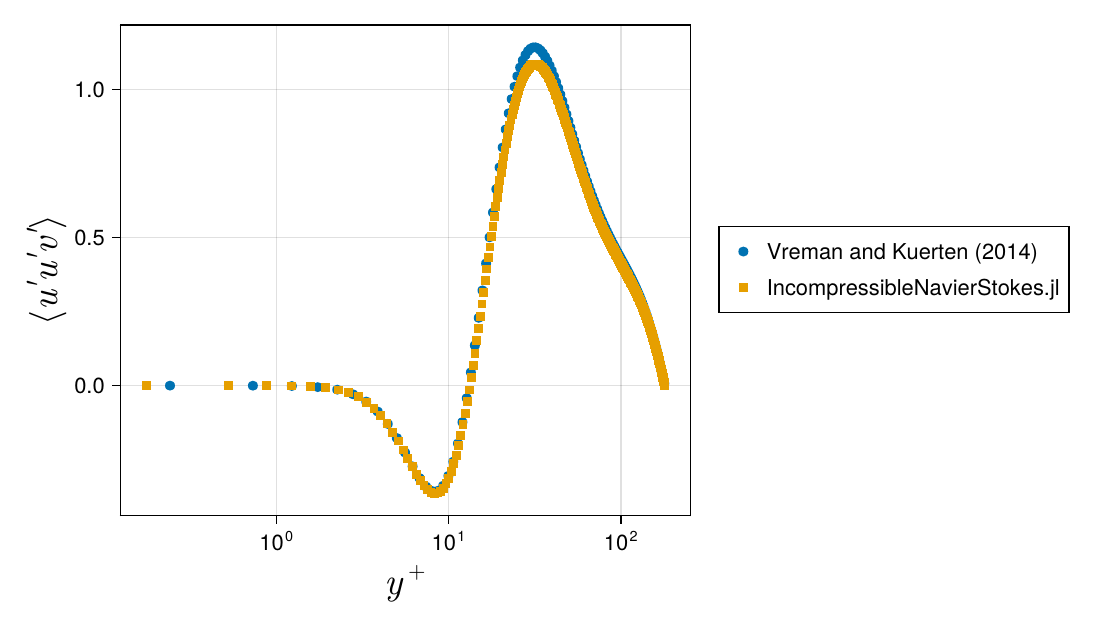}
    \caption{%
        Third order cross-moment of the fluctuations
        $\left\langle \overline{u' u' v'} \right\rangle$
        for DNS of turbulent channel flow
        with $512 \times 1024 \times 256$ grid points on a uniform
        grid, compared with reference data from Vreman and
        Kuerten~\cite{vremanComparisonDirectNumerical2014}.%
    }
    \label{S:fig:wallplot-dns-up2_vp}
\end{figure}

\begin{figure}
    \centering
    \includegraphics[width=\columnwidth]{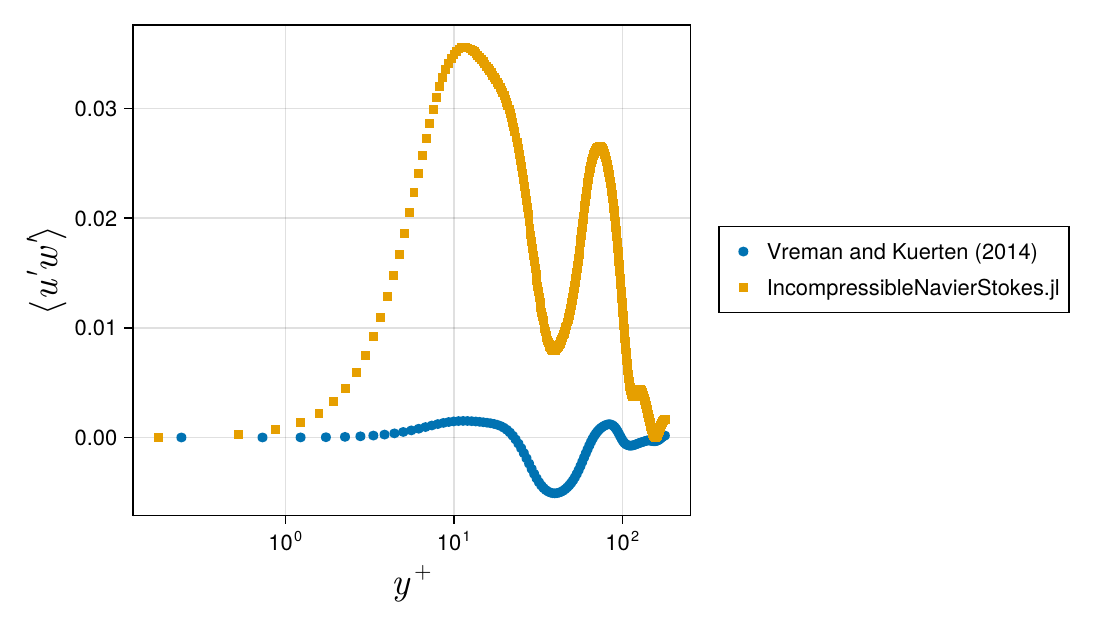}
    \caption{%
        Second order cross-moment of the fluctuations
        $\left\langle \overline{u' w'} \right\rangle$
        for DNS of turbulent channel flow
        with $512 \times 1024 \times 256$ grid points on a uniform
        grid, compared with reference data from Vreman and
        Kuerten~\cite{vremanComparisonDirectNumerical2014}.%
    }
    \label{S:fig:wallplot-dns-up_wp}
\end{figure}

Finally, \cref{S:fig:wallplot-dns-up2_vp,S:fig:wallplot-dns-up_wp} show the cross-moments
$\left\langle \overline{u' u' v'} \right\rangle$ and
$\left\langle \overline{u' w'} \right\rangle$.
Since the velocity components live on staggered face locations,
they are first interpolated to volume centers before computing the
cross-moments; the fluctuations $u' = u - \bar{u}$ are computed
before interpolation.
The cross-moment $\left\langle \overline{u' u' v'} \right\rangle$
is in good agreement with the reference data, apart from a small
underprediction at the peak consistent with the moments discussed above.
By contrast, $\left\langle \overline{u' w'} \right\rangle$ shows
no agreement with the reference.
This is expected: by the statistical symmetry of the channel in the
spanwise direction, $\left\langle \overline{u' w'} \right\rangle$
should vanish identically.
In practice, the finite averaging time leaves a residual of comparable
magnitude to the reference but with random sign, depending on the
initialization and floating-point perturbations.
The streamwise and wall-normal directions do not share this sensitivity,
as they are distinguished by the mean flow and the presence of the walls,
respectively.

Having validated the solver against DNS reference data, we now turn to
LES to demonstrate its ability to run standard eddy-viscosity closure models.

\subsection{LES with eddy-viscosity models}

We then perform LES of the same channel flow on a coarser grid with
$N_x = 128$, $N_y = 64$, $N_z = 64$ points.
The grid in the wall-normal direction ($y$) is non-uniform, using the
hyperbolic tangent stretching \eqref{S:eq:tanh-grid} to cluster points
near the walls. Since the grid is non-uniform, the spectral pressure
solver cannot be used; instead, the pressure Poisson equation is solved
using a sparse direct solver with an $L D L^T$ factorization,
implemented via the \texttt{cuDSS} library on the GPU.
This factorization requires substantial memory, which prevented us from
running the DNS at full resolution on a non-uniform grid: neither the
RTX~4090 (24\,GB) nor the H100 (94\,GB) had sufficient memory for the
$L D L^T$ decomposition at the $512 \times 256 \times 256$ resolution
of Vreman and Kuerten.
We employ five different eddy-viscosity models:
no-model (zero eddy-viscosity),
Smagorinsky~\cite{smagorinskyGeneralCirculationExperiments1963},
WALE~\cite{nicoudSubgridScaleStressModelling1999}, and
QR~\cite{verstappenWhenDoesEddy2011},
and Vreman's model~\cite{vremanEddyviscositySubgridscaleModel2004}.
The full model expressions are provided in \cref{S:sec:eddy-viscosity}.

We use the constants
$C_S = 0.1$ for Smagorinsky,
$C_W = 0.5$ for WALE,
$C_Q = \sqrt{3 / 2} / \pi$ for QR, and
$C_V = \sqrt{2.5 C_S^2}$ for Vreman's model.
The local filter width is defined based on the grid size with the common
definition $(\Delta x \Delta y \Delta z)^{1 / 3}$.
For Vreman's model, the individual grid spacings are incorporated separately.
We note that the choice of constants depend on the choice of filter width definition.
Verstappen argues that the filter width should also account for the interpolation scheme,
potentially doubling the filter width (and thus halving the ``optimal'' values
for the constants)~\cite{verstappenWhenDoesEddy2011}.

\begin{figure}
    \centering
    \includegraphics[width=\columnwidth]{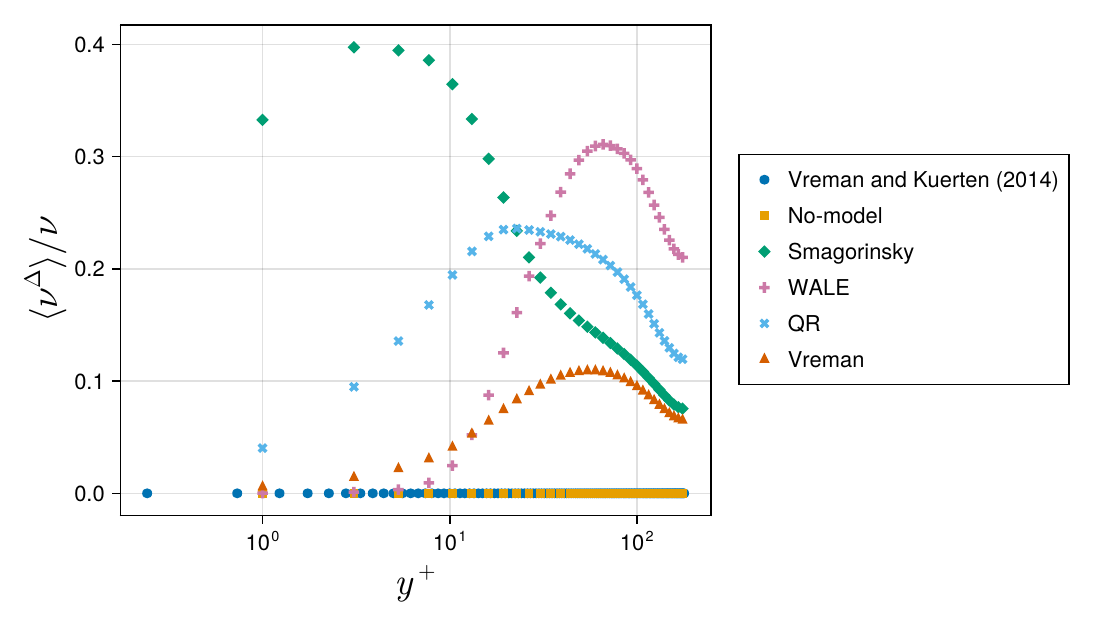}
    \caption{%
        Mean eddy-viscosity $\langle \nu^\Delta \rangle$
        normalized by the kinematic viscosity $\nu$
        for LES of turbulent channel flow
        with $128 \times 64 \times 64$ grid points on a non-uniform
        grid, compared with reference data from Vreman and
        Kuerten~\cite{vremanComparisonDirectNumerical2014}.%
        The DNS/no-model simulations do not produce any eddy-viscosity.%
    }
    \label{S:fig:wallplot-les-eddy-viscosity}
\end{figure}

\Cref{S:fig:wallplot-les-eddy-viscosity} shows the mean eddy-viscosity 
$\langle \bar{\nu}^\Delta \rangle / \nu$ normalized by the kinematic viscosity $\nu$.
The DNS and no-model simulations do not produce any eddy-viscosity (we set $\nu^\Delta = 0$).
The Smagorinsky viscosity does not decay to zero at the wall, which is a well-known issue with that model.
All the other models decay to zero at the wall, with the WALE model showing the fastest decay.
In the linear layer, Smagorinsky has the highest eddy-viscosity, followed by QR, Vreman, and WALE.
In the logarithmic region, WALE has the highest eddy-viscosity, followed by QR, Smagorinsky, and Vreman.
We emphasize that all the eddy-viscosities are scaled by the respective model constants. For different
choices of constants, the magnitude of the eddy-viscosities can change, but not the shape of the profiles.

\begin{figure}
    \centering
    \includegraphics[width=\columnwidth]{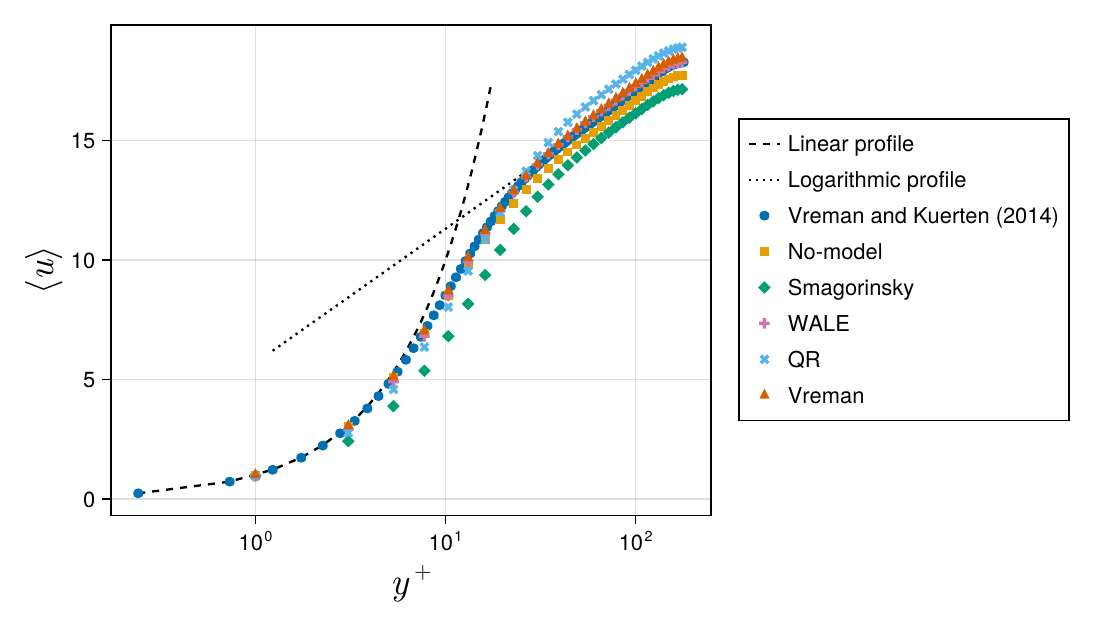}
    \caption{%
        Mean streamwise velocity profile $\langle \bar{u} \rangle$
        for LES of turbulent channel flow
        with $128 \times 64 \times 64$ grid points on a non-uniform
        grid, compared with reference data from Vreman and
        Kuerten~\cite{vremanComparisonDirectNumerical2014}.%
    }
    \label{S:fig:wallplot-les-ubar}
\end{figure}

\begin{figure}
    \centering
    \includegraphics[width=\columnwidth]{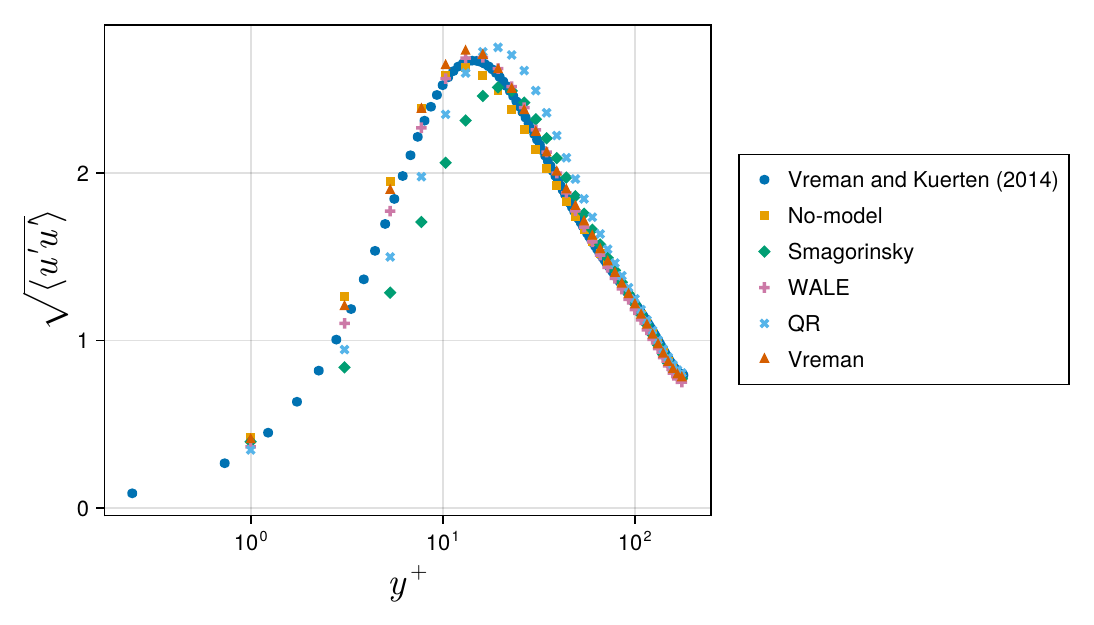}
    \caption{%
        Average squared velocity fluctuations in the $x$-direction
        for LES of turbulent channel flow
        with $128 \times 64 \times 64$ grid points on a non-uniform
        grid, compared with reference data from Vreman and
        Kuerten~\cite{vremanComparisonDirectNumerical2014}.%
    }
    \label{S:fig:wallplot-les-rms}
\end{figure}

\begin{figure}
    \centering
    \includegraphics[width=\columnwidth]{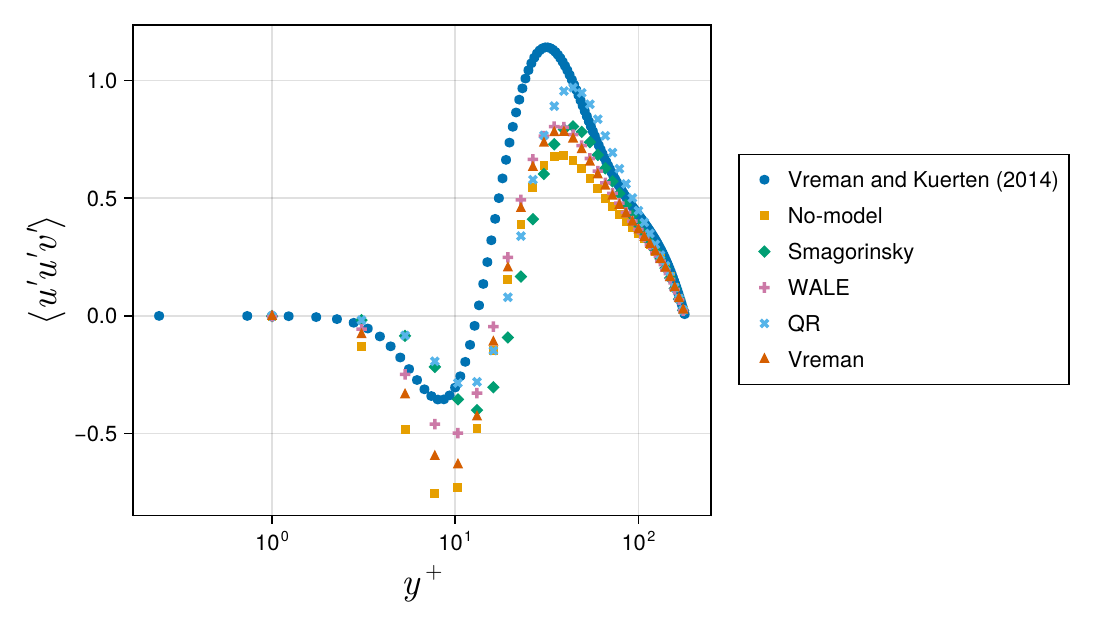}
    \caption{%
        Third order cross-moment of the fluctuations
        $\left\langle \overline{u' u' v'} \right\rangle$
        for LES of turbulent channel flow
        with $128 \times 64 \times 64$ grid points on a non-uniform
        grid, compared with reference data from Vreman and
        Kuerten~\cite{vremanComparisonDirectNumerical2014}.%
    }
    \label{S:fig:wallplot-les-up2_vp}
\end{figure}

\Cref{S:fig:wallplot-les-ubar,S:fig:wallplot-les-rms,S:fig:wallplot-les-up2_vp}
show the same statistics as in the DNS case, but for the LES.
For the average streamwise velocity profile $\langle \bar{u} \rangle$ in \cref{S:fig:wallplot-les-ubar},
Smagorinsky is too dissipative, resulting in an underprediction of
the velocity.
WALE and Vreman are closer to the reference than no-model.
The QR model is slightly below in the viscous layer and higher than the
reference in the logarithmic region.
For the root mean square streamwise velocity fluctuations
$\sqrt{\left\langle \overline{u' u'} \right\rangle}$ in \cref{S:fig:wallplot-les-rms},
WALE and Vreman perform better than no-model near the wall and near the channel
center, while no-model performs slightly better in the buffer layer.
This is consistent with the well-known tendency of eddy-viscosity models
to overdamp the turbulence production peak in the buffer layer.
Smagorinsky is again too dissipative, resulting in an underprediction of the
velocity fluctuations except for in the channel center.
For the cross-moment $\left\langle \overline{u' u' v'} \right\rangle$ in \cref{S:fig:wallplot-les-up2_vp},
the QR model performs best in the logarithmic region, while WALE is more
accurate near the wall.
All models perform better than no-model.

In conclusion, the solver correctly implements the eddy-viscosity models,
and the observed behavior is consistent with the
literature~\cite{nicoudSubgridScaleStressModelling1999}:
the Smagorinsky model is overly dissipative in wall-bounded flows due to
its non-vanishing eddy viscosity at the wall, while WALE, QR, and
Vreman's model all improve upon this by construction.
We note that the model constants and filter width definition were not
tuned for this case, and different choices can shift the
results~\cite{triasBuildingProperInvariants2015}.

\section{Hardware considerations and floating point arithmetic} \label{S:sec:hardware}

The choice of floating point precision has direct consequences for both the
accuracy and the computational cost of a simulation.

\subsection{Floating point formats and sources of error} \label{S:sec:float-representation}

Real numbers are represented on a computer using the IEEE~754 floating
point standard~\cite{ieeeStandard2019} as
$(-1)^s \times 2^{e - b} \times (1 + m)$, where $s$ is a sign bit,
$e$ is the exponent, $b$ is a format-dependent bias, and $m$ is the
significand (mantissa) with an implicit leading bit.
The formats most relevant for scientific computing are listed in
\cref{S:tab:float-formats}. Each format trades precision (significand
bits) and dynamic range (exponent bits) against storage size and
computational throughput. The machine epsilon
$\varepsilon = 2^{-(p-1)}$ bounds the relative error of a single
arithmetic operation.

\begin{table*}
    \centering
    \begin{tabular}{lcccccc}
        \toprule
        Format & Bits & Sign & Exponent & Significand
            & $\varepsilon$
            & Range \\
        \midrule
        Float64  & 64 & 1 & 11 & 52 & $\approx 10^{-16}$
            & $\approx 10^{\pm 308}$ \\
        Float32  & 32 & 1 & 8  & 23 & $\approx 10^{-7}$
            & $\approx 10^{\pm 38}$ \\
        Float16  & 16 & 1 & 5  & 10 & $\approx 10^{-3}$
            & $\approx 10^{\pm 4.5}$ \\
        BFloat16 & 16 & 1 & 8  & 7  & $\approx 10^{-2}$
            & $\approx 10^{\pm 38}$ \\
        \bottomrule
    \end{tabular}
    \caption{%
        Properties of common floating point formats. The machine epsilon
        $\varepsilon = 2^{-(p-1)}$, where $p$ is the number of
        significand bits plus the implicit leading bit, determines the
        smallest relative perturbation that is representable. The range
        indicates the approximate magnitude of the largest representable
        finite number.
    }
    \label{S:tab:float-formats}
\end{table*}

For DNS, the most relevant error mechanisms are
\emph{catastrophic cancellation}---the loss of significant digits when
subtracting nearly equal numbers, which directly affects finite difference
stencils---and \emph{overflow}, where results exceed the representable
range (Float16 overflows at approximately $6.5 \times 10^4$, which can be
exceeded by common physical quantities).
In long-running simulations, round-off errors accumulate as
$O(\sqrt{N_t} \, \varepsilon)$ for stable schemes. For chaotic flows,
even double-precision trajectories eventually diverge from the true
solution, but the practical question is whether the statistical
properties of the flow are affected by the precision.

\subsection{Precision trends in GPU hardware} \label{S:sec:precision-trends}

GPU manufacturers have progressively shifted transistor budgets toward
lower-precision arithmetic units, driven by deep learning
workloads~\cite{dongarraHardwareTrendsImpacting2024}.
In some architectures, Float64 throughput is only $1/32$ or $1/64$ of
Float32 throughput. Halving the number of bits per value also halves the
memory footprint and doubles the effective memory bandwidth, both of which
are critical for memory-bound stencil-based PDE solvers. However, the
reduced dynamic range of 16-bit formats can lead to overflow and loss of
significance, so scientists must carefully evaluate numerical stability at
reduced precision.
Reduced-precision strategies have been explored in climate
modelling~\cite{thornesUseScaledependentPrecision2017,kimpsonClimateChangeModelling2023},
where stochastic rounding can compensate for the loss of precision.

Mixed-precision strategies offer a middle ground: the bulk of the
computation is performed in reduced precision, while critical
accumulations (such as inner products in iterative solvers) use double
precision~\cite{highamMixedPrecisionStrategies2024}.
\texttt{INS.jl} supports choosing the floating point format, but does not
currently implement any mixed-precision strategies.

\subsection{Implications for DNS} \label{S:sec:implications-for-dns}

\begin{figure*}
    \centering
    \includegraphics[width=0.8\textwidth]{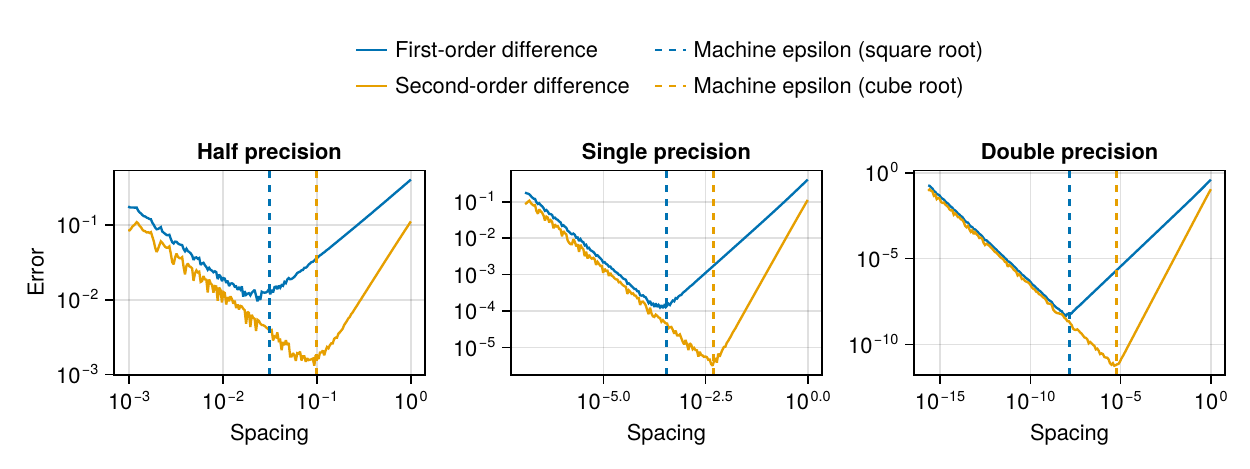}
    \caption{%
        Finite difference errors for the derivative of $\sin(x)$,
        evaluated in half (Float16), single (Float32), and double (Float64) precision.
        As the grid spacing $h$ decreases, the truncation error decreases
        until round-off error dominates, producing a characteristic
        V-shaped curve. The minimum error is reached at a grid spacing
        that depends on both the order of the finite difference and the
        machine epsilon of the floating point format. For second-order
        differences in single precision, the optimal spacing is comparable
        to typical DNS grid spacings, indicating that single precision is
        marginal at high resolutions.%
    }
    \label{S:fig:finite-difference-errors}
\end{figure*}

A finite difference approximation incurs both truncation error (decreasing
with grid spacing $h$) and round-off error (increasing with smaller $h$
due to catastrophic cancellation). Let
\begin{align}
    \partial^{h, 1}_x f(x) & \coloneq \frac{f(x + h) - f(x)}{h}, \\
    \partial^{h, 2}_x f(x) & \coloneq \frac{f(x + h) - f(x - h)}{2 h}
\end{align}
denote the first- and second-order finite difference approximations.
The error is minimized at $h \approx \varepsilon^{1/2}$ for the
first-order and $h \approx \varepsilon^{1/3}$ for the second-order
difference, as illustrated in \cref{S:fig:finite-difference-errors}.

For a DNS with $N = 1024$ points per direction in a domain of size
$2\pi$, the grid spacing is $h \approx 6 \cdot 10^{-3}$.
In single precision, the optimal spacing for second-order differences is
$h \approx 5 \cdot 10^{-3}$, which is comparable to the DNS spacing.
This suggests that single precision is marginal but feasible for
second-order methods at this resolution, consistent with the findings of
Karp et al.~\cite{karpEffectsLowerFloating2026}, who observed no
significant discrepancies between single- and double-precision results
across multiple turbulence solvers and test cases. At higher resolutions or
for higher-order methods, double precision becomes necessary.
Half precision is not viable for spatial discretization at practical DNS
resolutions, but may find application in components where high accuracy is
not required, such as neural network inference within a closure model.
In the closure modelling work of \cite{agdesteinDiscretizeFirstFilter2025},
neural networks are trained and evaluated in single precision,
but the tolerance of neural networks to quantization suggests
that half precision could be used for inference without significant loss of
closure model quality, while benefiting from the higher throughput of 16-bit
arithmetic units.

\subsection{Convergence study: the 2D Taylor-Green vortex} \label{S:sec:taylor-green}

The 2D Taylor-Green vortex is an exact solution to the incompressible
Navier-Stokes equations on the periodic domain $\Omega = [0, 2\pi]^2$.
Since the velocity and pressure fields are known in closed form,
this test case allows us to measure the spatial discretization error
directly.

We discretize the domain with $N \times N$ grid points for
$N \in \{4, 8, 16, \dots, 2048\}$ and compute the numerical solution at
$t = 1$ using a sufficiently small time step so that the temporal
discretization error is negligible compared to the spatial error.
The $L^2$ error between the numerical and analytical solutions
is then evaluated for three floating-point formats: Float64 (double
precision), Float32 (single precision), and Float16 (half precision).
We perform this convergence study on both a uniform grid and a
non-uniform grid with $\tanh$-stretching.

\begin{figure*}
    \centering
    \includegraphics[width=0.49\textwidth]{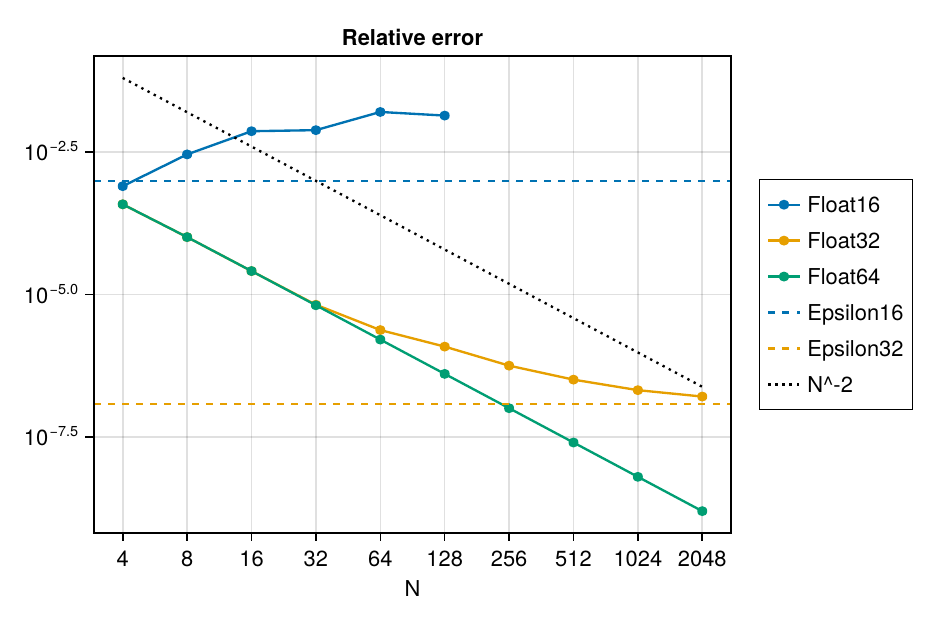}
    \includegraphics[width=0.49\textwidth]{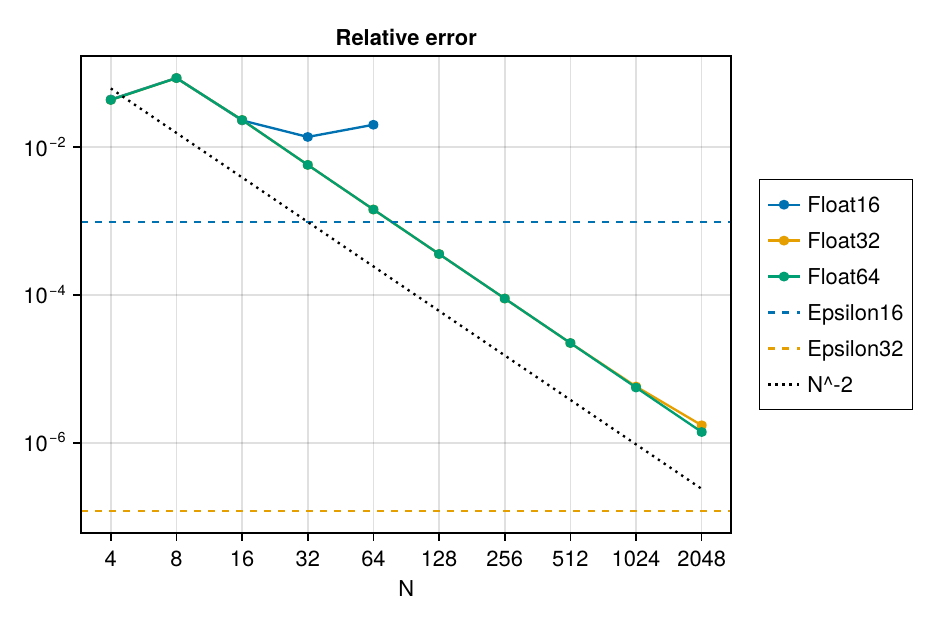}
    \caption{%
        Convergence of the numerical solution to the analytical solution for
        the velocity field of the 2D Taylor-Green vortex at $t = 1$.
        The error is measured in the
        $L^2$ norm and plotted against the number of grid points $N$ for different
        floating-point formats.
        The dotted line indicates second-order convergence.
        The dashed lines indicate the machine epsilon for single and half precision.
        Left: Uniform grid.
        Right: Non-uniform grid with $\tanh$-stretching.
    }
    \label{S:fig:taylor-green-convergence}
\end{figure*}

\paragraph{Uniform grid.}
On the uniform grid, the pressure Poisson equation is solved spectrally
using FFTs.
\Cref{S:fig:taylor-green-convergence} (left) shows that Float64 and
Float32 both exhibit clean second-order convergence, consistent with the
second-order finite-volume scheme. For Float32, the convergence rate
begins to deteriorate at fine resolutions and eventually plateaus as the
error reaches single-precision machine epsilon. Float16 does not display
a meaningful convergence trend and produces large errors across all
resolutions.

\paragraph{Non-uniform grid.}
On the non-uniform grid, the pressure is solved with a sparse direct
($LDL^T$) solver, which is only available for Float64 and Float32.
Since no direct solver implementation is available for Float16, the
Float16 case falls back to the conjugate gradient method in that case.
\Cref{S:fig:taylor-green-convergence} (right) shows that Float64 and
Float32 converge at the same rate, since the spatially varying grid
spacing leads to larger errors at a given $N$ compared to the uniform
case, so that the single-precision plateau is not yet reached for the
resolutions tested. Float16 tracks the other two formats at coarse
resolutions but becomes unstable as $N$ increases.

In summary, Float64 delivers reliable second-order convergence across all
tested resolutions. Float32 also converges at second order for coarse and
moderate resolutions, but the payoff of grid refinement diminishes at
finer grids as rounding errors accumulate and the effective convergence
rate drops below two. Float16, in its current form, is not viable for
this test case. One possible path toward making reduced-precision
arithmetic more effective is to scale the physical quantities so that they
remain well within the representable range of the floating-point format,
thereby reducing overflow, underflow, and loss of significance.

\section{Conclusion and outlook} \label{S:sec:conclusion}

We summarize the contributions of this work and outline directions for
future development.

\subsection{Summary}

This article presented \texttt{INS.jl}, an open-source
Julia package for DNS and LES of incompressible flows on staggered Cartesian
grids. The solver uses a second-order finite-volume spatial discretization
with explicit time integration and a pressure Poisson equation for
divergence-free projection, supporting both uniform and non-uniform grids
in two and three dimensions.

The software is designed around matrix-free, hardware-agnostic kernels that
are compiled for multi-threaded CPU or GPU from a single source
implementation. Each discrete operator is implemented in both a pure
(non-mutating) variant for automatic differentiation and a mutating variant
for memory-efficient forward simulation. Hand-written adjoint kernels
registered through \texttt{ChainRules.jl} enable reverse-mode differentiation through
the entire solver, making it possible to train neural network closure models
\emph{a posteriori} by backpropagating through the time integration.
Neural network closure models from \texttt{Lux.jl} are supported directly within the
solver, using the differentiable operator infrastructure for a-posteriori
training.

Three memory-optimization techniques---array reuse, low-storage
Runge--Kutta methods, and lazy computation of tensor components---allow
double-precision DNS at resolutions up to $840^3$ on a single datacenter GPU.

The code was validated against an analytical solution (the 2D Taylor--Green
vortex, confirming second-order convergence) and against reference DNS data
from Vreman and Kuerten~\cite{vremanComparisonDirectNumerical2014}
for turbulent channel flow at $\mathrm{Re}_\tau = 180$.
The mean velocity profile and RMS velocity fluctuations are in good agreement
with the reference. Higher-order moments (third and fourth) show small
discrepancies at the peaks, which we attribute to the second-order spatial
discretization.
LES of the same channel flow on a coarser non-uniform grid with four classical
eddy-viscosity models (Smagorinsky, WALE, QR, and Vreman) showed behavior
consistent with the literature: Smagorinsky is overly dissipative due to its
non-vanishing wall viscosity, while WALE, QR, and Vreman's model improve
upon this by construction, producing better agreement with the DNS reference
for both the mean velocity and velocity fluctuations.

We also investigated the effect of floating-point precision on accuracy:
Float64 delivers reliable second-order convergence, Float32 is effective at
moderate resolutions but offers diminishing returns at finer grids, and
Float16 is not yet viable in its current form.
Understanding these precision limits is important because modern GPU
hardware provides substantially higher throughput at reduced precision,
and neural network components within closure models may tolerate lower
precision for inference.

Software development practices---including version control, automated
testing, continuous integration, documentation generation, literate
programming for examples, release management with archival on Zenodo, and
reproducibility---ensure that the code remains maintainable and that
published results can be reproduced.
The code has been used by other
researchers~\cite{hoekstraReducedSubgridScale2026,
ivagnesNewDatadrivenEnergystable2026}, demonstrating its accessibility and
extensibility.

\subsection{Outlook}

The most immediate limitation of the current software is that it runs on a
single compute node. While the memory optimizations (\cref{S:sec:memory})
allow DNS at resolutions up to $840^3$ on a single datacenter GPU,
higher resolutions require distributing the computation across multiple GPUs
and multiple nodes. This requires domain decomposition with halo exchange---communication
of ghost-layer values at subdomain boundaries---between neighboring
subdomains, using CUDA-aware MPI for GPU-to-GPU
communication. The spectral pressure solver would need to be replaced with a
distributed variant based on pencil decomposition for the Fourier transforms.
The iterative conjugate gradient solver is more straightforward to distribute,
since it requires only local stencil operations and global reductions.

A second direction is to increase the spatial accuracy of the discretization.
The current second-order staggered scheme is well suited for DNS at moderate
resolutions, but higher-order stencils would reduce the number of grid points
needed to achieve a given accuracy. Verstappen and
Veldman~\cite{verstappenSymmetrypreservingDiscretizationTurbulent2003}
describe a fourth-order scheme that preserves the symmetry properties of the
second-order discretization by combining the standard discretization with a
coarser grid. Implementing this scheme would require extending the kernel
infrastructure to support wider stencils and multiple grid levels, while
preserving the energy-conservation properties that the current discretization
guarantees.

On the automatic differentiation side, the current approach uses \texttt{Zygote.jl}
for reverse-mode differentiation, which requires the pure (non-mutating)
operator variants and stores all intermediate states for the backward pass.
This limits the number of time steps that can be unrolled during training
before memory is exhausted. Checkpointing strategies, where only a subset
of intermediate states is stored and the rest are recomputed during the
backward pass, would allow training over longer time horizons without
proportionally increasing memory requirements. An alternative direction is
to use \texttt{Enzyme.jl}, which operates at a lower level than \texttt{Zygote.jl} and can
differentiate through mutating code directly, potentially avoiding the need
for separate pure and mutating operator variants altogether. Preliminary
support for \texttt{Enzyme.jl} is already included in the package as an extension
module.

The analysis in \cref{S:sec:hardware} showed that single precision is
marginal for second-order methods at typical DNS resolutions. Systematic
exploration of mixed-precision strategies---using reduced precision for
memory-bound kernels while retaining double precision for
accumulations---could improve performance on modern GPU hardware without
compromising accuracy.

Finally, the functional programming design of the solver makes it
straightforward to extract and reuse individual components in other
contexts. The discrete operators, grid generation routines, and initial
condition generators have already been reused in a pseudospectral solver
developed separately by the author. Maintaining this composability as the
package grows in scope will require careful attention to the software
architecture, ensuring that new features are added as modular components
rather than tightly coupled extensions.

\section*{Software and reproducibility statement}

The code used to generate the results of this paper is available at
\url{https://github.com/agdestein/IncompressibleNavierStokes.jl}.
It is released under the MIT license.

\section*{CRediT author statement}

\textbf{Syver Døving Agdestein}:
Conceptualization,
Formal analysis,
Methodology,
Software,
Validation,
Visualization,
Writing -- Original Draft

\textbf{Benjamin Sanderse}:
Funding acquisition,
Project administration,
Supervision,
Writing -- Review \& Editing

\section*{Declaration of Generative AI and AI-assisted technologies in the writing process}

During the preparation of this work the authors used Anthropic Claude Caude for
improving language and GitHub Copilot to propose wordings. After using these
tools/services, the authors reviewed and edited the content as needed and take
full responsibility for the content of the publication.

\section*{Declaration of competing interest}

The authors declare that they have no known competing financial interests or
personal relationships that could have appeared to influence the work reported
in this paper.

\section*{Acknowledgements}

This work is supported by the project ``Discretize first, reduce next'' (with
project number VI.Vidi.193.105) of the research programme NWO Talent Programme
Vidi  financed by the Dutch Research Council (NWO).
This work used the Dutch national e-infrastructure with the support of the
SURF Cooperative using grant no. EINF-15798.

\appendix

\section{Finite volume discretization} \label{P2:sec:discretization}

The integral form of the Navier-Stokes equations
serves as the starting point for developing a spatial discretization:
\begin{align}
    \int_{\partial \mathcal{O}}
    u \cdot n \, \mathrm{d} \Gamma & = 0, \label{P2:eq:mass_integral} \\
    \int_\mathcal{O}
    \left(\frac{\mathrm{d} u}{\mathrm{d} t} - f\right)
    \, \mathrm{d} \Omega
    & =
    \int_{\partial \mathcal{O}}
    \left( - u \otimes u - p I + \nu \nabla u \right) \cdot n \,
    \mathrm{d} \Gamma,
    \label{P2:eq:momentum_integral}
\end{align}
where $(u^1, \dots, u^d)$ is the velocity field, 
$p$ is the pressure field, $\nu$ is the viscosity,
$f$ is the external forcing,
$I$ is the identity matrix, and
$\mathcal{O} \subset \Omega$ is an arbitrary control volume with boundary
$\partial \mathcal{O}$, normal $n$, surface element $\mathrm{d} \Gamma$, and
volume size $|\mathcal{O}|$. We divide by the control volume sizes in the
integral form so that the system
\eqref{P2:eq:mass_integral}-\eqref{P2:eq:momentum_integral} has the same dimensions as
the continuous Navier-Stokes equations.

\subsection{Staggered grid configuration}

In this section, we describe a finite volume discretization of equations
\eqref{P2:eq:mass_integral}-\eqref{P2:eq:momentum_integral}. Before doing so, we
introduce our notation, which is such that the mathematical description of the
discretization closely matches the software implementation.

The $d$ spatial dimensions are indexed by $\alpha \in \{1, \dots, d\}$.
The $\alpha$-th unit vector is denoted $2 h_\alpha = (2 h_{\alpha \beta})_{\beta
= 1}^d$, where the (half) Kronecker symbol $h_{\alpha \beta}$ is $1 / 2$ if
$\alpha = \beta$ and $0$ otherwise. The Cartesian index $I = (I_1, \dots,
I_d)$ is used to avoid repeating terms and equations $d$ times, where
$I_\alpha$ is a scalar index (typically one of $i$, $j$, and $k$ in common
notation). This notation is dimension-agnostic, since we can write $u_I$ instead
of $u_{i j}$ in 2D or $u_{i j k}$ in 3D. In our Julia implementation of the
solver we use the same Cartesian notation (\verb|u[I]| instead of \verb|u[i, j]|
or \verb|u[i, j, k]|).

\begin{figure}
    \centering
    \includegraphics[width=\columnwidth]{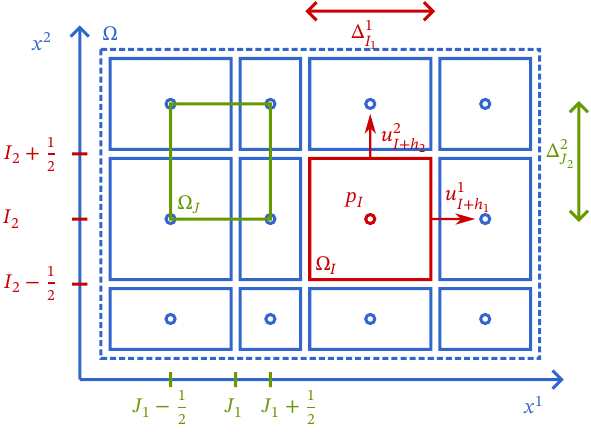}
    \caption{
        Finite volume discretization on a staggered grid. Note that the grid can
        be non-uniform, as long as each volume in a given column has the same
        width and each volume in a given row has the same height. Here, $I$ and
        $J$ are two arbitrary Cartesian indices, with $I \in \mathcal{I}$ in a
        volume center and $J \in \mathcal{I} + h_1 + h_2$ in a volume corner for
        illustrative purposes.
    }
    \label{P2:fig:finitevolumes}
\end{figure}

For the discretization, we use a staggered Cartesian grid as proposed by
Harlow and Welch~\cite{harlowNumericalCalculationTimeDependent1965}. Staggered
grids have excellent conservation
properties~\cite{lillyComputationalStabilityNumerical1965,
perotDiscreteConservationProperties2011}, and in particular their exact
divergence-freeness is essential. Consider a rectangular domain
$\Omega = \prod_{\alpha = 1}^d [a_\alpha, b_\alpha]$, where $a_\alpha <
b_\alpha$ are the domain boundaries and $\prod$ is a Cartesian product. Let
$\Omega = \bigcup_{I \in \mathcal{I}} \Omega_I$ be a partitioning of $\Omega$,
where $\mathcal{I} = \prod_{\alpha = 1}^d \{ \frac{1}{2}, 2 - \frac{1}{2},
\dots, N_\alpha - \frac{1}{2} \}$ are volume center indices, $N = (N_1, \dots,
N_d) \in \mathbb{N}^d$ are the number of volumes in each dimension,
$\Omega_I = \prod_{\alpha = 1}^d \Delta^\alpha_{I_\alpha}$ is a finite
volume, $\Gamma^\alpha_I = \Omega_{I - h_\alpha} \cap \Omega_{I + h_\alpha} =
\prod_{\beta \neq \alpha} \Delta^\beta_{I_\beta}$ is a volume face,
$\Delta^\alpha_i = \left[ x^\alpha_{i - \frac{1}{2}}, x^\alpha_{i + \frac{1}{2}}
\right]$ is a volume edge, $x^\alpha_0, \dots, x^\alpha_{N_\alpha}$ are volume
boundary coordinates, and $x^\alpha_i = \frac{1}{2} \left(x^\alpha_{i -
\frac{1}{2}} + x^\alpha_{i + \frac{1}{2}}\right)$ for $i \in \{ 1 / 2, \dots,
N_\alpha - 1 / 2\}$ are volume center coordinates. We also define the operator
$\delta_\alpha$ which maps a discrete scalar field $\varphi = (\varphi_I)_I$ to
\begin{equation} \label{P2:eq:discrete-derivative}
    (\delta_\alpha \varphi)_I = \frac{\varphi_{I + h_\alpha} - \varphi_{I -
    h_\alpha}}{| \Delta^\alpha_{I_\alpha} |}.
\end{equation}
It can be interpreted as a discrete equivalent of the continuous operator
$\frac{\partial}{\partial x^\alpha}$. All the above definitions are extended to
be valid in volume centers $I \in \mathcal{I}$, volume faces $I \in \mathcal{I}
+ h_\alpha$, or volume corners $I \in \mathcal{I} + \sum_{\alpha = 1}^d
h_\alpha$. The discretization is illustrated in figure~\ref{P2:fig:finitevolumes}.

\subsection{Equations for unknowns}

We now define the unknown degrees of freedom. The average pressure in
$\Omega_I$, $I \in \mathcal{I}$ is approximated by $p_I(t)$. The
average $\alpha$-velocity on the face $\Gamma^\alpha_I$, $I \in \mathcal{I} +
h_\alpha$ is approximated by $u^\alpha_I(t)$. Note that the pressure
$p$ and the $d$ velocity fields $u^\alpha$ are each defined at their own
canonical positions $x_I$ and $x_{I + h_\alpha}$ for $I \in \mathcal{I}$. This
is illustrated for a volume $I$ in figure~\ref{P2:fig:finitevolumes}. We now
derive the governing equations for these unknowns.

Using the pressure control volume $\mathcal{O} = \Omega_I$ with $I \in
\mathcal{I}$ in the integral constraint \eqref{P2:eq:mass_integral} and
approximating the face integrals with the mid-point quadrature rule
$\int_{\Gamma_I} u \, \mathrm{d} \Gamma \approx | \Gamma_I | u_I$ results in the
discrete divergence-free constraint
\begin{equation} \label{P2:eq:mass_discrete}
    \sum_{\alpha = 1}^d
    (\delta_\alpha u^\alpha)_I = 0.
\end{equation}
Note that dividing by the volume size yields a discrete equation resembling
the continuous one (since $| \Omega_I | = | \Gamma^\alpha_I | |
\Delta^\alpha_{I_\alpha} |$).

Similarly, choosing an $\alpha$-velocity control volume $\mathcal{O} =
\Omega_{I}$ with $I \in \mathcal{I} + h_\alpha$ in equation
\eqref{P2:eq:momentum_integral}, approximating the volume and face integrals using
the midpoint quadrature rule, and replacing the remaining spatial derivatives in
the diffusive term with a finite difference approximation yields the discrete
momentum equations
\begin{equation} \label{P2:eq:momentum_discrete}
\begin{split}
    \frac{\mathrm{d}}{\mathrm{d} t} u^\alpha_{I}
    & =
    - \sum_{\beta = 1}^d
    \left[\delta_\beta \left( \eta^\text{half}_\beta u^\alpha \eta^\text{lin}_\alpha u^\beta \right)\right]_{I} \\
    & + \nu \sum_{\beta = 1}^d (\delta_\beta \delta_\beta u^\alpha)_{I}
    + f^\alpha(x_{I})
    - (\delta_\alpha p)_{I}.
\end{split}
\end{equation}
where we have assumed that $f$ is constant in time for simplicity.
The outer discrete derivative in $(\delta_\beta \delta_\beta u^\alpha)_{I}$ is
required at position $I$, which means that the inner derivative is evaluated
at $(\delta_\beta u^\alpha)_{I + h_\beta}$ and $(\delta_\beta u^\alpha)_{I -
h_\beta}$, thus requiring $u^\alpha_{I - 2 h_\beta}$, $u^\alpha_{I}$, and
$u^\alpha_{I + 2 h_\beta}$, which are all at their canonical positions. The two
velocity components $\eta^\text{half}_\beta u^\alpha$
and $\eta^\text{lin}_\alpha u^\beta$
in the convective term 
are required at
positions $I - h_\beta$ and $I + h_\beta$, which are outside the canonical
positions. Their values at the required positions are obtained using averaging with
weights $1 / 2$ for the $\alpha$-component and linear interpolation for the
$\beta$-component. This preserves the skew-symmetry of the convection operator,
ensuring energy conservation in the convective term~\cite{vangastelenEnergyconservingNeuralNetwork2024}.

\section{Discrete operator stencils} \label{S:sec:stencils}

This appendix lists the explicit stencil expressions for the discrete
operators introduced in \cref{S:sec:spatial}, along with their adjoint
(pullback) kernels used for reverse-mode automatic differentiation
(\cref{S:sec:differentiability}).
The adjoint kernels are tested against numerical finite-difference gradients to
verify their correctness.

All operators are expressed in terms of the discrete derivative
$\delta_\alpha$ defined in \eqref{P2:eq:discrete-derivative}. Below, $I$
denotes a multi-index at a grid point, $h_\alpha = e_\alpha / 2$ is a
half-index shift, $|\Delta^\alpha_i|$ is the width of volume $i$ in the
$\alpha$-direction, and $\bar{\varphi}$, $\bar{u}$, $\bar{p}$ denote
adjoint (cotangent) variables. All pullbacks are derived from the general
definition \eqref{S:eq:pullback}.
Except for convection, all operators
are linear, and their pullbacks reduce to the transpose of the operator.

\paragraph{Divergence}
The divergence maps the velocity field $u$ to a scalar field $\varphi$
defined at the pressure points. The kernel and its pullback are
\begin{equation} \label{S:eq:div-kernel}
    \varphi_I = \sum_{\alpha = 1}^d
    \frac{u^\alpha_{I + h_\alpha} - u^\alpha_{I - h_\alpha}}
    {|\Delta^\alpha_{I_\alpha}|},
    \quad
    \bar{u}^\alpha_{I} =
    \frac{\bar{\varphi}_{I - h_\alpha}}
    {|\Delta^\alpha_{I_\alpha - 1/2}|}
    - \frac{\bar{\varphi}_{I + h_\alpha}}
    {|\Delta^\alpha_{I_\alpha + 1/2}|}.
\end{equation}

\paragraph{Pressure gradient}
The pressure gradient maps a scalar field $p$ to a vector field $\varphi$
defined at the velocity points. The kernel and its pullback are
\begin{equation} \label{S:eq:presgrad-kernel}
    \varphi^\alpha_I
    = \frac{p_{I + h_\alpha} - p_{I - h_\alpha}}
    {|\Delta^\alpha_{I_\alpha}|},
    \quad
    \bar{p}_{I} = \sum_{\alpha = 1}^d
    \left(
    \frac{\bar{\varphi}^\alpha_{I - h_\alpha}}
    {|\Delta^\alpha_{I_\alpha - 1/2}|}
    - \frac{\bar{\varphi}^\alpha_{I + h_\alpha}}
    {|\Delta^\alpha_{I_\alpha + 1/2}|}
    \right).
\end{equation}

\paragraph{Diffusion}
The diffusion maps a vector field $u$ to a vector field $\varphi$ defined
at the velocity points. The kernel and its pullback are
\begin{equation} \label{S:eq:diff-kernel}
    \varphi^\alpha_I =
    \nu
    \sum_{\beta = 1}^d
    \frac{1}{|\Delta^\beta_{I_\beta}|}
    \left(
    \frac{u^\alpha_{I + 2 h_\beta} - u^\alpha_{I}}
    {|\Delta^\beta_{I_\beta + 1/2}|}
    - \frac{u^\alpha_{I} - u^\alpha_{I - 2 h_\beta}}
    {|\Delta^\beta_{I_\beta - 1/2}|}
    \right)
\end{equation}
and
\begin{equation} \label{S:eq:diff-pullback}
\begin{split}
    \bar{u}^\alpha_{I} =
    \nu
    \sum_{\beta = 1}^d
    & \Bigg(
    \frac{\bar{\varphi}^\alpha_{I - 2 h_\beta}}
    {|\Delta^\beta_{I_\beta - 1}| \, |\Delta^\beta_{I_\beta - 1/2}|} \\
    - & \frac{1}{|\Delta^\beta_{I_\beta}|}
    \left(
    \frac{\bar{\varphi}^\alpha_{I}}
    {|\Delta^\beta_{I_\beta + 1/2}|}
    + \frac{\bar{\varphi}^\alpha_{I}}
    {|\Delta^\beta_{I_\beta - 1/2}|}
    \right) \\
    + & \frac{\bar{\varphi}^\alpha_{I + 2 h_\beta}}
    {|\Delta^\beta_{I_\beta + 1}| \, |\Delta^\beta_{I_\beta + 1/2}|}
    \Bigg).
\end{split}
\end{equation}

\paragraph{Convection}
The convective term maps a vector field $u$ to a vector field $\varphi$
defined at the velocity points.
The product $u^\alpha u^\beta$ at off-diagonal positions requires
interpolation of velocity components to the appropriate face locations.
The interpolation weights
\begin{equation} \label{S:eq:interp-weights}
\begin{split}
    A^{\alpha \beta +}_{i}
    & = \frac{|\Delta^\beta_{i - \delta_{\alpha \beta} + 1}|}
    {|\Delta^\beta_{i - \delta_{\alpha \beta}}|
    + |\Delta^\beta_{i - \delta_{\alpha \beta} + 1}|}, \\
    A^{\alpha \beta -}_{i}
    & = \frac{|\Delta^\beta_{i + \delta_{\alpha \beta} - 1}|}
    {|\Delta^\beta_{i + \delta_{\alpha \beta} - 1}|
    + |\Delta^\beta_{i + \delta_{\alpha \beta}}|}
\end{split}
\end{equation}
map to the right ($+$) and left ($-$) in the $\beta$-direction,
respectively, and are chosen to preserve
skew-symmetry~\cite{verstappenSymmetrypreservingDiscretizationTurbulent2003}.
The convective kernel is
\begin{equation} \label{S:eq:conv-kernel}
\begin{split}
    \varphi^\alpha_I = -\sum_{\beta = 1}^d
    \frac{1}{|\Delta^\beta_{I_\beta}|}
    \bigg( \\
    \left(
        A^{\alpha \beta +}_{I_\beta} u^\alpha_{I} +
        A^{\alpha \beta -}_{I_\beta + 1} u^\alpha_{I + 2 h_\beta}
    \right) \\
    \left(
        A^{\beta \alpha +}_{I_\alpha} u^\beta_{I} +
        A^{\beta \alpha -}_{I_\alpha + 1} u^\beta_{I + 2 h_\alpha}
    \right) \\
    - \left(
        A^{\alpha \beta +}_{I_\beta - 1} u^\alpha_{I - 2 h_\beta} +
        A^{\alpha \beta -}_{I_\beta} u^\alpha_{I}
    \right) \\
    \left(
        A^{\beta \alpha +}_{I_\alpha - \delta_{\alpha \beta}} u^\beta_{I - 2 h_\beta} +
        A^{\beta \alpha -}_{I_\alpha - \delta_{\alpha \beta} + 1}
        u^\beta_{I - 2 h_\beta + 2 h_\alpha}
    \right)
    \bigg).
\end{split}
\end{equation}
Since the convective kernel is nonlinear, the pullback contains components
from the primal input $u$. The input $u$ appears eight times in
\eqref{S:eq:conv-kernel}, and their contributions to the pullback are
\begingroup
\small
\begin{equation} \label{S:eq:conv-pullback}
    \begin{split}
    & \bar{u}^\alpha_{I} = \\
    & -\bar{\varphi}^\alpha_{I}
    \sum_{\beta = 1}^d
    \frac{1}{|\Delta^\beta_{I_\beta}|}
    A^{\alpha \beta +}_{I_\beta}
    \left(
        A^{\beta \alpha +}_{I_\alpha} u^\beta_{I} +
        A^{\beta \alpha -}_{I_\alpha + 1} u^\beta_{I + 2 h_\alpha}
    \right) \\
    & -\bar{\varphi}^\alpha_{I - 2 h_\beta}
    \sum_{\beta = 1}^d
    \frac{1}{|\Delta^\beta_{I_\beta - 1}|}
    A^{\alpha \beta -}_{I_\beta}
    \left(
        A^{\beta \alpha +}_{I_\alpha - \delta_{\alpha \beta}} u^\beta_{I - 2 h_\beta} +
        A^{\beta \alpha -}_{I_\alpha - \delta_{\alpha \beta} + 1}
        u^\beta_{I - 2 h_\beta + 2 h_\alpha}
    \right) \\
    & -\sum_{\beta = 1}^d
    \bar{\varphi}^\beta_{I}
    \frac{1}{|\Delta^\alpha_{I_\alpha}|}
    \left(
        A^{\beta \alpha +}_{I_\alpha} u^\beta_{I} +
        A^{\beta \alpha -}_{I_\alpha + 1} u^\beta_{I + 2 h_\alpha}
    \right)
    A^{\alpha \beta +}_{I_\beta} \\
    & -\sum_{\beta = 1}^d
    \bar{\varphi}^\beta_{I - 2 h_\beta}
    \frac{1}{|\Delta^\alpha_{I_\alpha - \delta_{\alpha \beta}}|}
    \left(
        A^{\beta \alpha +}_{I_\alpha - \delta_{\alpha \beta}} u^\beta_{I - 2 h_\beta} +
        A^{\beta \alpha -}_{I_\alpha - \delta_{\alpha \beta} + 1}
        u^\beta_{I - 2 h_\beta + 2 h_\alpha}
    \right)
    A^{\alpha \beta -}_{I_\beta} \\
    & +\bar{\varphi}^\alpha_{I + 2 h_\beta}
    \sum_{\beta = 1}^d
    \frac{1}{|\Delta^\beta_{I_\beta + 1}|}
    A^{\alpha \beta +}_{I_\beta + 1}
    \left(
        A^{\beta \alpha +}_{I_\alpha + \delta_{\alpha \beta}} u^\beta_{I + 2 h_\beta} +
        A^{\beta \alpha -}_{I_\alpha + \delta_{\alpha \beta} + 1}
        u^\beta_{I + 2 h_\beta + 2 h_\alpha}
    \right) \\
    & +\bar{\varphi}^\alpha_{I}
    \sum_{\beta = 1}^d
    \frac{1}{|\Delta^\beta_{I_\beta}|}
    A^{\alpha \beta -}_{I_\beta + 1}
    \left(
        A^{\beta \alpha +}_{I_\alpha} u^\beta_{I} +
        A^{\beta \alpha -}_{I_\alpha + 1} u^\beta_{I + 2 h_\alpha}
    \right) \\
    & +\sum_{\beta = 1}^d
    \bar{\varphi}^\beta_{I + 2 h_\alpha}
    \frac{1}{|\Delta^\alpha_{I_\alpha + 1}|}
    \left(
        A^{\beta \alpha +}_{I_\alpha + 1} u^\beta_{I + 2 h_\alpha} +
        A^{\beta \alpha -}_{I_\alpha + 2} u^\beta_{I + 4 h_\alpha}
    \right)
    A^{\alpha \beta +}_{I_\beta + \delta_{\alpha \beta}} \\
    & +\sum_{\beta = 1}^d
    \bar{\varphi}^\beta_{I - 2 h_\beta + 2 h_\alpha}
    \frac{1}{|\Delta^\alpha_{I_\alpha - \delta_{\alpha \beta} + 1}|} \\
    & \left(
        A^{\beta \alpha +}_{I_\alpha - \delta_{\alpha \beta} + 1}
        u^\beta_{I - 2 h_\beta + 2 h_\alpha} +
        A^{\beta \alpha -}_{I_\alpha - \delta_{\alpha \beta} + 2}
        u^\beta_{I - 2 h_\beta + 4 h_\alpha}
    \right)
    A^{\alpha \beta -}_{I_\beta + \delta_{\alpha \beta}}.
    \end{split}
\end{equation}
\endgroup

\section{Eddy-viscosity models} \label{S:sec:eddy-viscosity}

This appendix lists the eddy-viscosity closure models implemented in the
software, following the unified invariant-based framework of Trias et
al.~\cite{triasBuildingProperInvariants2015}.

\subsection{General framework}

All eddy-viscosity models express the sub-grid stress as
\begin{equation} \label{S:eq:eddy-viscosity}
    m^\alpha(u) = -\sum_{\beta = 1}^d \partial_\beta (2 \nu_t S_{\alpha \beta}),
\end{equation}
where
$S_{\alpha \beta} = (\partial_\beta u^\alpha + \partial_\alpha u^\beta) / 2$
is the resolved strain-rate tensor and $\nu_t \geq 0$ is the eddy viscosity.

On the staggered grid, $S_{\alpha \beta}$ is naturally defined at the
staggered tensor points: the diagonal components $\bar{S}_{\alpha \alpha}$
live at the pressure points, while the off-diagonal components
$\bar{S}_{\alpha \beta}$ ($\alpha \neq \beta$) live at the edge midpoints.
The eddy viscosity $\nu_t$ is computed at the pressure points from the
velocity gradient $A_{\alpha \beta} = \delta_\beta u^\alpha$, which is
interpolated from the staggered tensor points to the pressure points.
To evaluate \eqref{S:eq:eddy-viscosity}, $\nu_t$ must then be interpolated
back to the staggered tensor points where $S_{\alpha \beta}$ is
defined. For the diagonal components, $\nu_t$ is already at the correct
location. For the off-diagonal components, $\nu_t$ is interpolated using a
four-point linear average of the surrounding pressure points.

\subsection{Invariants of the velocity gradient tensor}

All models below are expressed in terms of the velocity gradient tensor
$A = \nabla u$ at the pressure points and its symmetric and
anti-symmetric parts,
\begin{equation}
    S = \frac{A + A^\mathsf{T}}{2}, \qquad
    W = \frac{A - A^\mathsf{T}}{2}.
\end{equation}
Following Trias et al.~\cite{triasBuildingProperInvariants2015}, we
define the invariants
\begin{equation} \label{S:eq:invariants}
\begin{split}
    Q_A = -\frac{\operatorname{tr}(A^2)}{2}, \quad
    Q_S = -\frac{\operatorname{tr}(S^2)}{2}, \quad
    Q_W = -\frac{\operatorname{tr}(W^2)}{2}, \\
    R_S = \frac{\operatorname{tr}(S^3)}{3}, \quad
    R_A = \frac{\operatorname{tr}(A^3)}{3},
\end{split}
\end{equation}
and the mixed invariant
\begin{equation} \label{S:eq:V2}
    V^2 = 4 \left( \operatorname{tr}(S^2 W^2) - 2 Q_S Q_W \right).
\end{equation}
The five independent invariants of $A$ are $Q_S$, $Q_W$, $R_S$,
$R_A$, and $V^2$. Note that $Q_A = Q_S + Q_W$ and that $-2 Q_S =
S_{\alpha \beta} S_{\alpha \beta} \geq 0$, so $Q_S \leq 0$.

\subsection{Implemented models}

All models take the form $\nu_t = (C \Delta)^2 \, \mathcal{D}$, where
$C$ is a dimensionless model constant and $\Delta$ is the filter width.
The quantity $\mathcal{D}$ has dimensions of inverse time and differs
between models. The proposed values of $C$ are taken from the original
references or from Trias et
al.~\cite{triasBuildingProperInvariants2015}.

\paragraph{Smagorinsky}
The Smagorinsky
model~\cite{smagorinskyGeneralCirculationExperiments1963,lillyRepresentationSmallscaleTurbulence1966}
defines the eddy viscosity as
\begin{equation} \label{S:eq:smagorinsky}
    \nu_t = (C_\text{S} \Delta)^2 \sqrt{2 S_{\alpha \beta} S_{\alpha \beta}}
    = (C_s \Delta)^2 \sqrt{-4 Q_S},
\end{equation}
with proposed constant $C_\text{S} = 0.17$. This is the simplest and most widely
used eddy-viscosity model. Its main limitation is that it does not vanish
at solid walls or in laminar flow regions, requiring ad-hoc damping
functions or dynamic procedures.

\paragraph{Vreman}
The Vreman
model~\cite{vremanEddyviscositySubgridscaleModel2004}
defines
\begin{equation} \label{S:eq:vreman}
    \nu_t = C_\text{Vr}^2
    \sqrt{\frac{B_\beta}{\operatorname{tr}(A A^T) / 2}},
\end{equation}
where $\beta_{i j} = \Delta^2_m A_{i m} A_{j m}$ accounts for non-isotropic
grid and
$B_\beta = \beta_{1 1} \beta_{2 2} - \beta_{1 2}^2 +
\beta_{1 1} \beta_{3 3} - \beta_{1 3}^2 +
\beta_{2 2} \beta_{3 3} - \beta_{2 3}^2$.
The proposed constant is $C_\text{Vr} = \sqrt{2.5 C_\text{S}^2}$.
The denominator $\operatorname{tr}(A A^\mathsf{T}) / 2 \geq 0$ vanishes only
when $A = 0$, so the model is well-defined wherever
the velocity gradient is nonzero.
The viscosity is set to zero when $A = 0$.

\paragraph{QR (Verstappen)}
Verstappen's minimum-dissipation QR
model~\cite{verstappenWhenDoesEddy2011,triasBuildingProperInvariants2015}
defines
\begin{equation} \label{S:eq:qr}
    \nu_t = -(C_\text{Ve} \Delta)^2 \frac{|R_S|}{Q_S},
\end{equation}
with proposed constant $C_\text{Ve} = \sqrt{3 / 2} / \pi \approx 0.390$~\cite{verstappenWhenDoesEddy2011}
or $C = 0.527$~\cite{triasBuildingProperInvariants2015}.
This model vanishes in two-component flows (where $R_S = 0$), which is
the correct physical behavior near solid walls and in regions of
axisymmetric strain.

\paragraph{WALE}
The wall-adapting local eddy-viscosity (WALE)
model~\cite{nicoudSubgridScaleStressModelling1999,triasBuildingProperInvariants2015}
defines
\begin{equation} \label{S:eq:wale}
    \nu_t = (C_\text{W} \Delta)^2
    \frac{\left( \mathcal{S}^d_{i j} \mathcal{S}^d_{i j} \right)^{3/2}}
    {\left(S_{i j} S_{i j}\right)^{5/2}
    + \left( \mathcal{S}^d_{i j} \mathcal{S}^d_{i j} \right)^{5/4}},
\end{equation}
where $\mathcal{S}^d$ is the traceless symmetric part of $A^2$: $\mathcal{S}^d \coloneq (A A + A^T A^T) / 2 - (Q_A / 3) I$
with proposed constant $C_\text{W} = \sqrt{2.5 C_\text{S}}$, where $C_\text{S}$ is the Smagorinsky constant.
The WALE model is designed to produce the correct $y^+$ scaling of the
eddy viscosity near solid walls, without requiring explicit damping
functions.

\paragraph{$\sigma$-model (Nicoud)}
The $\sigma$-model~\cite{nicoudUsingSingularValues2011,triasBuildingProperInvariants2015}
is defined in terms of the singular values $\sigma_1 \geq \sigma_2 \geq
\sigma_3 \geq 0$ of the velocity gradient tensor $A$:
\begin{equation} \label{S:eq:sigma}
    \nu_t = (C_\sigma \Delta)^2
    \frac{\sigma_3 (\sigma_1 - \sigma_2)(\sigma_2 - \sigma_3)}
    {\sigma_1^2},
\end{equation}
with proposed constant $C_\sigma = 1.35$.
This model vanishes for two-component and axisymmetric flows, and
automatically produces the correct near-wall scaling without damping
functions.

\paragraph{S3PQR (Trias)}
The S3PQR model~\cite{triasBuildingProperInvariants2015} is a family of
models parameterized by an exponent $p$:
\begin{equation} \label{S:eq:s3pqr}
    \nu_t = (C_\text{Tr} \Delta)^2
    P_{AA}^{\, p} \, Q_{AA}^{-(p+1)} \, R_{AA}^{(p + 5/2)/3},
\end{equation}
where
\begin{equation}
\begin{split}
    P_{AA} & = 2(Q_W - Q_S) = \operatorname{tr}(A A^\mathsf{T}), \\
    Q_{AA} & = V^2 + Q_A^2, \\
    R_{AA} & = R_A^2.
\end{split}
\end{equation}
The choice of $p$ recovers or approximates several known models.
Trias et al.\ propose $p \in \{-5/2, \, -1, \, 0\}$ as representative
members of the family.

\end{document}